   \def\obs#1{{\bf (*** #1 ***)} }

 \def\obs#1{}     
\documentclass[twoside,letterpaper,draft,11pt]{amsart}

\usepackage{amsmath,amsthm,latexsym,xspace,amscd,amssymb,xypic}               

\usepackage{color}

\title[Partial Galois cohomology and  related homomorphisms]{Partial Galois cohomology and related homomorphisms (expanded version)}

\author[M.\ Dokuchaev]{M. Dokuchaev}
\address{Instituto de
Matem\'atica e Estat\'\i stica,
Universidade de S\~ao Paulo, Rua do Mat\~ao, 1010,
05508-090 S\~ao Paulo, SP, Brasil}
\email{dokucha@ime.usp.br}

\author[A.\ Paques]{A. Paques}
\address{Instituto de
Matem\'atica  Universidade  Federal do Rio Grande do Sul  Avenida Bento Gon\c calves 9500,
91509-900 Porto Alegre, RS, Brasil}
\email{paques@mat.ufrgs.br}

\author[H.\ Pinedo ]{H. Pinedo }
\address{Escuela de
Matem\'aticas, Universidad Industrial de Santander, Cra 27 calle 9, Edificio Camilo Torres, Bucaramanga, Colombia.}
\email{hpinedot@uis.edu.co}

\thanks{{The first author was partially supported by FAPESP and CNPq of Brazil. The second author was  partially supported by FAPESP  of Brazil. The third author was supported by FAPESP of Brazil.}
}

\newtheorem{teo}{Theorem}[section]
\newtheorem{defi}[teo]{Definition}
\newtheorem{lema}[teo]{Lemma}

\newtheorem{prop}[teo]{Proposition}
\newtheorem{claim}[teo]{Claim}
\newtheorem{exe}[teo]{Example}
\newtheorem{remark}[teo]{Remark}

\newcommand{\La}{\Lambda}

\newcommand{\0}{\theta}
\newcommand{\Om}{\Omega}

\newcommand{\ep}{\epsilon}

\newcommand{\de}{\delta}

\newcommand{\cua}{\hfill${\blacksquare}$}
\newcommand{\omgh}{\omega_{g,h}}
\newcommand{\af}{\alpha}
\newcommand{\afg}{\alpha_g}
\newcommand{\afh}{\alpha_h}

\newcommand{\afgh}{\alpha_{gh}}

\newcommand{\lb}{\lambda}

\newcommand{\te}{\theta}

\newcommand{\p}{{\bf Proof. }}

\newcommand{\om}{\omega}

\newcommand{\mfp}{{\mathfrak p}}

\newcommand{\U}{{\mathcal U}}
\newcommand{\df}{\displaystyle\frac}
\newcommand{\ot}{{\otimes}}

\newcommand{\m}{{}^{-1}}
\newcommand{\mt}{\mapsto}
\def\f{\varphi}
\def\e{\varepsilon}

\def\ndv{\ {\mid \kern -0.7 em {\scriptstyle \not}} \ \ }

\def\nd{\ {\mid \kern -0.4 em {\scriptstyle \not}} \ \ }

\newcommand{\N}{{\mathbb N}}

\newcommand{\vm}{\vspace{.05cm}}
\newcommand{\vu}{\vspace{.1cm}}
\newcommand{\vd}{\vspace{.2cm}}
\newcommand{\vt}{\vspace{.3cm}}

\begin{document}

\date{\today}

\begin{abstract}  For a partial Galois extension of commutative rings we give a seven terms sequence, which is an analogue of the Chase-Harrison-Rosenberg   sequence.   
\end{abstract}
\maketitle
\begin{section}{Introduction}

  The concept of a Galois extension of commutative rings was introduced in the same paper by M. Auslander and O. Goldman \cite{AG}, in which they  laid the foundations for  separable extensions of commutative rings and defined the Brauer group of a commutative ring.  Later, in  \cite{CHR},  S. U. Chase, D.K. Harrison and  A. Rosenberg developed  Galois  theory of commutative rings by giving several equivalent definitions of a Galois extension, establishing a Galois correspondence, and  specifying, to the case of a Galois extension, the Amitsur cohomology seven terms exact sequence,  given by S. U. Chase and  A. Rosenberg in \cite{CR}.  The Chase-Harrison-Rosenberg   sequence  can be viewed as a common generalization of the two most fundamental facts from Galois co\-ho\-mo\-lo\-gy of fields: the Hilbert's Theorem 90 and the isomorphism of  the relative Brauer group with the second cohomology group of the Galois group. Since then much attention have been payed to the  sequence and its parts subject to more constructive proofs, generalizations and analogues in various  contexts.
  
  Our point of view is to replace global actions by partial ones. The latter are becoming  an object of intensive research and have their origins  in the theory of operator algebras, where they, together with the corresponding  crossed products and partial representations,  form the essential ingredients of a new and successful method to study   $C^*$-algebras generated by partial isometrics, initiated by R. Exel  in \cite{E-1}, \cite{E-2},   \cite{E0},  \cite{E2} and  \cite{E1}.  The first algebraic results on these new concepts, established in \cite{E1}, \cite{DEP}, \cite{St1}, \cite{St2}, \cite{KL} and \cite{DE}, and the development of a  Galois theory of partial actions in~\cite{DFP}, stimulated a growing algebraic activity around  partial actions  (see the surveys~\cite{D3} and~\cite{F2}). In particular,   partial Galois theoretic results  have been obtained in \cite{BP}, \cite{CaenDGr},  \cite{CaenJan},  \cite{FrP},  \cite{KuoSzeto}, \cite{PRSantA}, and  applications of partial actions were found to graded algebras in~\cite{DE} and \cite{DES},  to tiling semigroups in~\cite{KL1},   to Hecke algebras in~\cite{E3},   to automata theory in \cite{DNZh},   to restriction semigroups in~\cite{CornGould}  and \cite{Kud} and to  Leavitt path algebras in~\cite{GR}. In addition, the interpretation  of the famous R. Thompson's groups as partial action groups on finite binary words permitted J.-C. Birget  to study algorithmic problems for them  \cite{Birget}. Amongst the recent advances we mention a remarkable application of the theory of partial actions to paradoxical decompositions and to algebras related to separated graphs \cite{AraE1}, its efficient use  in the study of the Carlsen-Matsumoto $C^*$-algebras associated to arbitrary subshifts~\cite{DE2} and of the  Steinberg algebras~\cite{BeuGon2}, as well as the proof of an algebraic version of the Effros-Hahn conjecture on the ideals in partial crossed products~\cite{DE3}.

 The general notion of a continuous twisted partial action of a locally compact group on a $C^*$-algebra introduced in \cite{E0}, and adapted to the abstract ring theoretic context in \cite{DES}, contains multipliers which satisfy a sort of $2$-cocycle identity, and it was natural to ask   
   what kind of cohomology theory would suite it. The answer was  given in \cite{DK}, where the partial cohomology of groups were introduced and studied together with their relation to cohomology of inverse semigroups, showing also that it  fits nicely the theory of partial projective group representations developed in \cite{DN}, \cite{DN2} and \cite{DoNoPi}.  Note that partial group cohomology turned out to be  useful to study ideals of global  reduced $C^*$-crossed products  \cite{KennedySchafhauser}.
   
   Having at hand  partial Galois theory and partial group cohomology we may ask now what would be the  analogue of the Chase-Harrison-Rosenberg  exact  sequence in the context of a partial Galois extension of commutative rings. The purpose of the present paper  is to  answer this question.   The additional new ingredients include 
a partial action of the Galois group $G$ on the disjoint union of the Picard groups of all direct summands of  $R$ (see Section~\ref{phi1phi2phi3}), as well as partial representations of $G$ (see Section~\ref{phi6}).

   In Section~\ref{prelim} we recall for reader's convenience some facts used in the paper, whereas the homomorphisms of the sequence are given in  Sections~\ref{phi1phi2phi3}, ~\ref{phi4phi5}, ~\ref{phi6}.
   
   Our proofs are constructive and the partial case is essentially more laborious than the classical one, so that in this   article  we only build up the homomorphisms,  and the proof of the exactness of the sequence will be given in a forthcoming   paper.

Throughout this   work  the word ring  means an associative ring with an identity element. For any ring $R$ by an $R$-module we mean a left unital  $R$-module.
If $R$ is commutative,  we shall consider an $R$-module $M$ as a central $R$-$R$-bimodule, i.e, an $R$-$R$-bimodule $M$
with $mr = rm$ for all $m\in M$ and $r\in R.$   
 We write that $M$ is a {\it  f.g.p. $R$-module} if $M$ is a (left) projective and finitely generated $R$-module,      and by a {\it  faithfully projective
 $R$-module} we mean  a faithful, f.g.p.  $R$-module. For a monoid (or a ring)
 $T,$ the group  of its units (i.e, invertible elements) is denoted by $\U(T).$  In all what follows, unless otherwise stated, $R$ will denote a commutative ring and unadorned $\ot$ will mean $\ot_R.$

\end{section}
 
\section{Preliminaries}\label{prelim}

In this section we give some definitions and results which will be used  in  the paper. All modules over commutative rings are considered as central
bimodules.  
\vu

\subsection{The Brauer group of a commutative ring}

  Recall that an $R$-algebra $A$ is called \emph{separable} if $A$ is a projective module over its enveloping algebra $A^e=A\ot A^{\rm{op}},$ where $A^{\rm{op}}$ denotes the opposite algebra of $A.$
If $A$ is faithful as an $R$-module we can
identify $R$ with $R1_A$, and if, in addition, its center $C(A)$ is equal to $R$ we say that $A$ is \emph{central}. Moreover, $A$
is called an \emph{Azumaya} $R$-algebra if $A$ is central and separable. Equivalently, $A$ is a faithfully projective $R$-module and
$A\ot A^{\rm{op}} \simeq {\rm End}_R (A)$ as $R$-algebras,  (see \cite[Theorem 2.1(c)]{AG}).

 \vu

 In \cite{AG} the following  equivalence relation was defined on the class of all Azumaya $R$-algebras:

\vu

 \noindent $A\sim B$ if there exist faithfully projective $R$-modules $P$ and $Q$ such that
$$A \otimes {\rm End}_R(P)\cong B \otimes {\rm End}_R(Q), $$ as $R$-algebras.

\vu

 Let $[A]$ denote the equivalence class containing $A$, and $B(R)$ the set of all such equi\-va\-len\-ce classes.
  Then $B(R)$ has a natural structure of a multiplicative abelian group, whose multiplication is induced by the tensor product of $R$-algebras,
 that is, $$[A][B]=[A \otimes B],\, \ \text{for all}\,\  [A],[B]\in B(R).$$ Its identity element is $[R]$ and $[A]^{-1}=[A^{\rm{op}}]$,
 for all $[A]\in B(R)$.  This group is called the {\it  Brauer group of R.}

 \vu

According to \cite{AG}, for any commutative $R$-algebra $S,$  the map from $B(R)$ to $B(S)$, given by $[A]\mapsto [A\otimes S]$, is a
well defined group homomorphism. Its kernel is denoted by  $B(S/R)$ and called the {\it relative Brauer  group of S over $R$}.
If $A$ is an Azumaya $R$-algebra  whose equivalence class in $B(R)$ belongs to  $B(S/R),$  we say that $A$ is {\it split} by $S,$ or
that $S$ is a {\it splitting ring} for $A.$

\vu

For any nonempty subset $X$ of a ring $B$ and any subring $V$ of $B,$ we denote by $C_V(X)=\{y\in V\,|\, xy= yx\, {\rm \,\,for\, \,all}\, \,x\in X\}$
the {\it centralizer} of $X$ in $V.$ In particular, if $X=V,$ then $C_V(V)$ is the {\it center} $C(V)$  of $V.$ It is known that a
commutative $R$-subalgebra $B$ of an $R$-algebra $A$ is a {\it maximal commutative subalgebra} if and  only if $C_A(B)=B.$

\vu

\subsection{Partial cohomology of groups}  Let $G$ be a group.
A   {\it unital twisted  partial action}  of  $G$ on $R$ is a triple $$\af=(\{D_g\}_{g\in G}, \{\afg\}_{g\in G}, \{\omgh\}_{(g,h)\in G\times G}),$$
such that for every $g\in G,$ $D_g$ is an ideal of $R$ generated by a non-necessarily non-zero  idempotent $1_g,$ $\afg\colon D_{g\m}\to D_g$
is a ring isomorphism, for each pair $(g,h)\in G\times G,$ $\omgh\in \U(D_gD_{gh})$ and  for all $g,h,l\in G$ the following statements are satisfied:

\smallskip

\noindent (i) $D_1=R$ and $\af_1$ is the identity map of $R,$

\noindent (ii) $\afg(D_{g\m}D_h)=D_gD_{gh},$

\noindent (iii) $\afg \circ \afh(t)=\omgh\afgh(t)\omgh^{-1}$ for any $t\in D_{h\m} D_{(gh)\m},$

\noindent (iv) $\om_{1,g}=\om_{g,1}=1_g$  and

\noindent (v) $\afg(1_{g\m}\om_{h,l})\om_{g,hl}=\omgh\om_{gh,l}.$

\smallskip

Using (v) one obtains
\begin{equation}\label{afom}\afg(\om_{g\m,g})=\om_{g,g\m},\,\, \,\,\text{for any}\,\,g\in G.\end{equation}

In \cite{DES} the authors defined twisted partial actions of groups on algebras in a more general setting in which the $D_g$'s are
not necessarily unital rings. In all what follows we shall use only unital twisted partial actions.

\vu

If  $(\{D_g\}_{g\in G}, \{\afg\}_{g\in G}, \{\om_{g,l}\}_{(g,l)\in G\times G})$ is a twisted partial action of $G$
on $R,$ the family of partial isomorphisms $(D_g,\afg)_{g\in G}$ forms a  {\it partial action}\footnote {In the sense defined in \cite{DE}.}
which we denote by $\af.$  Then, the family $\om=\{\omgh\}_{(g,h)\in G\times G}$ is called a {\it twisting} of $\af,$ and the
above twisted partial action will be denoted by $(\af, \om).$

\vu

If $R$, in particular, is a multiplicative monoid,  then  one obtains  from the above definition the concept of a unital twisted partial
action of a group on a monoid.
\smallskip

We recall from \cite{DK} the following.
\begin{defi}\label{defn-cochain}
Let  $T$ be a commutative ring or  monoid,  $n\in\N, n>0$ and $\af=(T_g, \af_g)_{g\in G}$ a unital partial action of $G$ on T. An $n$-cochain
of $G$ with  values in $T$ is a function $f:G^n\to T,$ such that $f(g_1,\dots,g_n)\in\U(T1_{g_1}1_{g_1g_2}\dots1_{g_1g_2\dots g_n}).$
A $0$-cochain is an element of $\U(T)$.
\end{defi}

\vu

Denote the set of all $n$-cochains by $C^n(G,\af, T)$.  This set is an abelian group via the point-wise multiplication.
Its identity is the map  $(g_1,\dots , g_n)\mapsto  1_{g_1}1_{g_1g_2}\dots1_{g_1g_2\dots g_n} $  and the inverse of $f\in C^n(G,\af,T)$ is
$f^{-1}(g_1,\dots,g_n)=f(g_1,\dots,g_n)^{-1}$, where $f(g_1,\dots,g_n)^{-1}$ is the inverse of $f(g_1,\dots,g_n)$
in $T1_{g_1}1_{g_1g_2}\dots1_{g_1g_2\dots g_n},$ for all $g_1,\dots, g_n\in G$.

\vu

\begin{defi}{\rm (}{\bf The coboundary homomorphism}{\rm)}
Given $n\in\N, n>0,$ $f\in C^n(G,\af, T)$ and $g_1,\dots,g_{n+1}\in G$, set
\begin{align}\label{pcob}
(\delta^nf)(g_1,\dots,g_{n+1})=&\,\af_{g_1}\left(f(g_2,\dots,g_{n+1})1_{g\m_1}\right)
\prod_{i=1}^nf(g_1,\dots , g_ig_{i+1}, \dots,g_{n+1})^{(-1)^i}\notag\\
&f(g_1,\dots,g_n)^{(-1)^{n+1}}.
\end{align}
Here, the inverse elements are taken in the corresponding ideals. If $n=0$ and $t$ is an invertible element of $T$, we set
$(\delta^0t)(g)=\af_g(t1_{g^{-1}})t^{-1},$ for all $g\in G$.
\end{defi}

\vm

\begin{prop}\cite[Proposition 1.5]{DK}\label{pcobh}
		$\delta^n$ is a group homomorphism from $C^n(G,\af, T)$ to $C^{n+1}(G,\af, T)$ such that
	$
		(\delta^{n+1}\delta^nf)(g_1,g_2,\dots, g_{n+2})=1_{g_1}1_{g_1g_2}\dots1_{g_1g_2\dots g_{n+2}},
	$
for any $n\in \N,$ $f\in C^n(G,\af,T)$ and $g_1,g_2,\dots, g_{n+2} \in G.$ 
	\end{prop}

\vm

\begin{defi}\label{defn-cohomology}
For	$n\in\N,$ we define the groups $Z^n(G,\af,T)=\ker{\delta^n}$ of  partial $n$-cocycles, $B^n(G,\af,T)$$={\rm im}\,{\delta^{n-1}}$
of partial $n$-co\-boun\-da\-ries, and $H^n(G,\af,T)=\df{\ker{\delta^n}}{{\rm im}\,{\delta^{n-1}}}$ of partial $n$-co\-ho\-mo\-lo\-gies of
$G$ with values in $T$, $n\ge 1.$ For $n=0$ we define $H^0(G,\af,T)=Z^0(G,\af,T)=\ker{\delta^0}$.
	\end{defi}

\begin{exe}
	\begin{align*}
		H^0(G,\af,T)&=Z^0(G,\af,T)=\{t\in{\mathcal U}(T)\mid{\af}_g(t 1_{g^{-1}}) = t 1_g,  \forall g \in G \},\\
		B^1(G,\af,T)&=\{ f\in C^1(G,\af,T) \mid f(g) = {\af}_g(t 1_{g^{-1}})  t^{-1}, \, \mbox{for some}\,\, t \in  {\mathcal U}(T)  \} .
	\end{align*}

We have	$ ( \delta^{1} f) (g,h) = \af _g ( f(h) 1_{g^{-1}} ) f(gh)^{-1}f(g)$
for $f \in  C^1 (G,\af,T),$ so that
	\begin{align*}
		Z^1(G,\af, T) &=\{f\in C^1 (G,\af,T)\mid f(gh)1_g=f(g) \, \af _g ( f(h) 1_{g^{-1}} ), \forall g,h\in G \},\end{align*} moreover
		$B^2 (G,\af,T)$ is the group \begin{align*}
&\left\{ w \in  C^2 (G,\af, T) \mid \exists  f \in  C^1 (G,\af, T),\text{with}\,\,w(g,h) = \af _g ( f(h) 1_{g^{-1}} )f(g) f(gh)\m
\right\}.
	\end{align*}
For $n=2$ we obtain
	$$
		( \delta^{2} w) (g,h,l) = \af _g ( w(h,l) 1_{g^{-1}} )  \; w(gh,l)^{-1} \; w(g,hl) \; w(g,h)^{-1},
	$$
with $w \in C^2 (G,\af,T),$ and $Z^2(G,\af,T)$ is
	\begin{align*}
		&\{w\in C^2 (G,\af, T) \mid \af _g ( w(h,l) 1_{g^{-1}} ) \; w(g,hl) = w(gh,l)  \, w(g,h),\,\,\forall g,h,l \in G \}.
	\end{align*}
Hence,  the elements of $Z^2(G,\af,T)$ are exactly the twistings for $\af.$ \cua \end{exe}

\vm

Two cocycles $f,f'\in Z^n(G,\af,T)$  are called {\it cohomologous}  if they differ by an $n$-coboundary.

\begin{remark}\label{normalized} Notice that a $1$-cocycle is always normalized, i.e. $f(1) = 1_T.$ Indeed, taking $g=h=1$ in the $1$-cocycle equality we immediately see that $f(1)= f(1)^2,$ so  $f(1)$ must be $1_T,$ as $f(1) \in \U (T).$   
\end{remark}

\vu

\subsection {Partial Galois extensions}

Let  $G$ be a finite group and $\alpha=(D_g, \afg)_{g\in G}$ a unital (non twisted) partial action of $G$ on $R$. The subring of {\it invariants} of
$R$ under $\af$ was introduced in \cite{DFP} as
\begin{equation} \label{inva}R^\af=\{r\in R\,|\, \afg(r1_{g\m})=r1_g \,\,{\rm for \,\, all} \, g\in G\}.\end{equation}  Notice that $R^\af=H^0(G,\af,T).$

The ring extension $R\supseteq R^\af$ is called an {\it $\af$-partial Galois extension} if for some $m\in \N$ there exists a
subset $\{x_i, y_i\,|\, 1\le i\le m\}$ of $R$ such that $\displaystyle\sum_{i=1}^{m}x_i\af_g(y_i 1_{g\m})=\delta_{1,g}, \,\,g\in G.$

 As in \cite{DFP}, we call the set $\{x_i, y_i\,|\, 1\le i\le m\}$ a {\it partial Galois coordinate system} of  $R\supseteq R^\af.$

\vu

The {\it trace map} ${\rm tr}_{R/R^\af}:R\to R^\af$ is given by  $x\mapsto\sum_{g\in G}\afg(x1_{g\m}).$ By \cite[Remark 3.4]{DFP}
there exists $c\in R$ such that ${\rm tr}_{R/R^\af}(c)=1 ,$ provided that the extension $R\supseteq R^\af$ is $\af$-partial Galois 
	 (here we write $1=1_R = 1_{R^\af}$).

\vu

The \emph{partial skew group ring} $R\star_\af G$ is defined as the set of all formal sums $\sum_{g\in G} r_g\de_g$, $r_g\in D_g,$ with
the usual addition and the multiplication determined by the rule $$(r_g\de_g)(r'_h\de_h)= r_g\af_g(r'_h1_{g\m})\de_{gh}.$$

It is shown in \cite[Theorem 4.1]{DFP} that  $R\supseteq R^\af$ is a partial Galois extension if and only if $R$ is a f.g.p.
$R^\af$-module and  the map

\begin{equation}\label{jota}
j\colon R\star_\af G \to {\rm End}_{R^\af}(R),\,\,\,\,\, j\left (\sum_{g\in G} r_g\de_g\right)(r)=\sum_{g\in G} r_g\afg(r1_{g\m})
\end{equation}
is an $R^\af$-algebra and an $R$-module isomorphism.

\vu

\section{On generalizations of the Picard group}\label{Picard}
To construct our version of the seven term  sequence, we need some generalizations of the concept of the Picard group of a commutative ring.
First, we recall  the next. 
\begin{defi} The abelian group of all $R$-isomorphism classes of  f.g.p. $R$-modules of rank 1, with binary operation given by $[P][Q]=[P\ot Q]$
is denoted by ${\bf Pic}(R).$ The identity in ${\bf Pic}(R)$ is  $[R]$ and the inverse of $[P]$ in ${\bf Pic}(R)$ is $[P^*],$ where
$M^*={\rm Hom}_R(M,R)$ for any $R$-module $M.$  \end{defi}

Recall that if $P$ is a faithfully projective   $R$-module, then $[P] \in  {\bf Pic}(R)$ exactly when the map $R\to {\rm End}_R (P),$
given by $r \mapsto m_r$,  where $ m_r(p)=r p$ for all $p\in P$,  is an isomorphism of rings  (see \cite[Lemma I.5.1]{DI}).
We also recall the next.

\begin{prop}\label{ht}\cite[Hom-Tensor Relation I.2.4]{DI}  Let  $A$ and $B$ be $R$-algebras. Let $M$ be a f.g.p. $A$-module and
$N$ be a f.g.p. $B$-module. Then for any $A$-module $M'$ and any $B$-module $N'$, the map
 $$\psi\colon {\rm Hom}_A(M,M')\ot {\rm Hom}_B(N,N')\to {\rm Hom}_{(A\ot B)}(M\ot N,M'\ot N' ),$$
 induced by $(f\ot  g)(m\ot n)=f(m)\ot g(n),$ for all $m\in M,\, n\in N$, is an $R$-module isomorphism. If $M=M'$ and $N=N',$
 then $\psi$ is an $R$-algebra isomorphism.\cua\end{prop}

Let  $\Lambda$ be a unital commutative $R$-algebra.  We give the following.

\begin{defi} A  $\La$-$\La$-bimodule $P$ is called  $R$-partially invertible if $P$ is   central as an $R$-$R$-bimodule, and
\begin{itemize}
\item P is a  left f.g.p. $\La$-module,
\item There is an $R$-algebra epimorphism $\Lambda^{\rm{op}} \to  {\rm End}_\Lambda(P).$
\end{itemize}

\vm

Let $$[P]=\{ M\,|\, M\, \text{is a}\,\, \La\text{-}\La\text{-bimodule and} \,\, M\cong P\, \,\text{as} \,\,\La\text{-}\La\text{-bimodules}  \}.$$
We denote by  ${\bf PicS}_R(\La)$  the set of the isomorphism classes $[P]$ of $R$-partially invertible $\La$-$\La$-bimodules. Finally, we set
${\bf PicS}_R (R):={\bf PicS} (R).$
\end{defi}

\begin{prop}\label{ccom}   The product $[P][Q]=[P\ot_\La Q]$ endows   ${\bf PicS}_{R}(\La)$ with the structure of a  semigroup. 
\end{prop}
\p We  shall show that  $[P\ot_\La Q]\in {\bf PicS}_{R}(\La),$ for any $[P], [Q]\in {\bf PicS}_{R}(\La).$  Notice that $P\ot_\La Q$ is a
left f.g.p. $\La$-module. Indeed, there are free f.g left  $\La$-modules $F_1,F_2$ and $\La$-modules $M_1,M_2$ such that
$P\oplus M_1=F_1,\, Q\oplus M_2=F_2.$ Now consider $M_1$ and $F_1$ as central $\La$-$\La$-bimodules, then by tensoring the two previous
relations we see that there exists a left $\La$-module $M$ such that $(P\ot_\La Q)\oplus M=F_1\ot_\La F_2,$ and the assertion follows.

\vu

By assumption there are $R$-algebra epimorphisms $\xi_1\colon \Lambda^{\rm{op}} \to  {\rm End}_\Lambda(P)$ and  $\xi_2\colon \Lambda^{\rm{op}} \to  {\rm End}_\Lambda(Q).$
It follows from Proposition \ref{ht} that $\xi_1\ot \xi_2  \colon \Lambda^{\rm{op}}\ot_\Lambda
\Lambda^{\rm{\rm{op}}}\to {\rm End}_{\La}(P\ot_\La Q)$ is an  $R$-algebra epimorphism. Since $\Lambda^{\rm{\rm{op}}}\ni \lb \mt \lb\ot_\Lambda 1_\Lambda \in
\Lambda^{\rm{op}}\ot_\Lambda \Lambda^{\rm{op}} $ is an $R$-algebra isomorphism, we conclude that  $ \Lambda^{\rm{op}} \ni \lambda\mapsto
\xi_1(\lb)\ot \xi_2(1_\Lambda )\in {\rm End}_{\La}(P\ot_\La Q)$ is an $R$-algebra epimorphism.\cua

\vu

 Throughout the paper, by ${\rm Spec}(R)$ we mean, as usual, the set of all prime ideals of $R$.

\begin{defi} We say that a f.g.p. central $R$-$R$-bimodule  $P$ has rank less than or equal to one, if for any  $\mfp\in {\rm Spec}(R)$
one has $P_\mfp=0$ or $P_\mfp\cong R_\mfp$ as $R_\mfp$-modules. In this case we write ${\rm rk}_{R}(P)\le 1.$
\end{defi}
The following result characterizes ${\bf PicS} (R).$

\begin{prop}\label{carpic} We have
\begin{equation*}{\bf PicS} (R)=\{[E]\,|\, E\, \text{is a f.g.p. central  $R$-$R$-bimodule and}\,\, {\rm rk}_{R}(E)\le 1\}.
\end{equation*}
\end{prop}
\p Let $E$ be a  f.g.p. central  $R$-$R$-bimodule such that ${\rm rk}_{R}(E)\le 1,$ and consider the map $$m_R \colon R \ni r
\mapsto m_r\in {\rm End}_R(E),\,\, m_r(x)=rx,\, r\in R, x\in E.$$ Via localization it is easy to show that $m_R$ is an $R$-algebra epimorphism.

\vu

Conversely if $[E]\in {\bf PicS} (R),$ then $E$ is a f.g.p. central  $R$-$R$-bimodule and there is an   $R$-algebra epimorphism
 $R\to  {\rm End}_R(E).$ Thus, for any $\mfp\in {\rm Spec}(R)$ there is an $R$-algebra epimorphism
$$R_{\mfp}\to {\rm End}_{R_{\mfp}}  (R^{n_\mfp}_{\mfp})\simeq M_{{n_\mfp}} (R_{\mfp}),$$ where $n_\mfp={\rm rk}_{R_\mfp}(E_\mfp),$ which gives $n_\mfp\le 1$.
\cua
\medskip

\begin{remark}\label{unit}  Notice that  $\U( {\bf PicS}(R))={\bf Pic}(R).$ Indeed, the inclusion   $\U( {\bf PicS}(R))\supseteq{\bf Pic}(R)$ is trivial.
On the other hand, for $[E] \in \U( {\bf PicS}(R))$ there exists $[P] \in  {\bf PicS}(R)$ with $E\otimes P \cong R,$ so that  $E_\mfp \otimes
P_\mfp  \cong R_\mfp $ for each prime $\mfp $ in $R.$ Then  $E_\mfp \neq 0$ and we see  by Proposition~\ref{carpic} that
${\rm rk}_{R_\mfp}\, (E_\mfp)=1$ for each prime $\mfp ,$ and  thus  $[E] \in  {\bf Pic}(R).$
\end{remark}

Given  an inverse semigroup $S,$ we denote the inverse of $s\in S$ by  $s^*.$  We proceed with the following fact.
\begin{prop}\label{picinv}The set  ${\bf PicS} (R)$ with the binary operation induced by the tensor product is a commutative inverse
monoid with 0. Moreover $[E^*]=[E]^*,$ for all $[E]\in {\bf PicS} (R).$
\end{prop}
\p It follows from Proposition \ref{ccom} and Proposition \ref{carpic}   that  ${\bf PicS} (R)$ is  a commutative monoid with $0.$

Take $[M]\in {\bf PicS} (R).$ By Proposition \ref{ht} we obtain $(M^*)_\mfp\cong (M_\mfp)^*={\rm Hom}_{R_\mfp}(M_\mfp, R_\mfp),$ for
all $\mfp \in {\rm Spec}(R),$ and hence  $[M^*]\in {\bf PicS} (R)$ thanks to Proposition \ref{carpic}.

Now we  prove that $[M][M^*][M]=[M]$ and $[M^*][M][M^*]=[M^*].$ Recall that $M \otimes M^* \cong {\rm End}_{R}(M),$
since $M$ is a f.g.p. $R$-module (see \cite[Lemma I.3.2 (a)]{DI}), and  we get $[M][M^*][M]=$$[{\rm End}_{R}(M)][M].$
There is an $R$-module homomorphism $$\kappa\colon {\rm End}_{R}(M) \ot  M\ni f\ot m\mapsto f(m)\in M,$$ and via localization we
will prove that $\kappa$ is an isomorphism.  Indeed, take $\mfp\in {\rm Spec}(R)$ then there are two cases to consider.

\emph{Case 1: $M_\mfp=0$.} In this case $\kappa_\mfp \colon 0\to 0$ is clearly an $R_\mfp$-module isomorphism.

\emph{Case 2: $M_\mfp\cong R_\mfp.$} Here we have   $\kappa_\mfp \colon R_p\ot_{R_\mfp} R_p\ni r'_\mfp\ot_{R_\mfp} r_\mfp\to
r'_\mfp r_\mfp\in R_\mfp$ is  an $R_\mfp$-module isomorphism.

From this we conclude that $[M][M^*][M]=[M],$ for all $[M]\in {\bf PicS} (R).$ Finally since $M$ is a f.g.p.
$R$-module, there is an $R$-module isomorphism $M\cong (M^*)^* $ (see \cite[Theorem V.4.1]{MC}), and
consequently    $[M^*][M][M^*]=[M^*][(M^*)^*][M^*]=[M^*].$\cua
\bigskip

By Proposition \ref{picinv} and Clifford's Theorem (see for instance \cite{CP}),   ${\bf PicS} (R)$ is a semilattice of
abelian groups. In particular, $${\bf PicS} (R)=\bigcup\limits_{\zeta \in F(R)}{\bf PicS}_\zeta (R),$$
where $F(R)$ is a semilattice isomorphic to the semilattice of the idempotents of ${\bf PicS} (R).$ Therefore,  to describe
${\bf PicS} (R)$ we need to know its idempotents.
\smallskip

   We recall that  given an inverse semigroup $S,$ its  idempotents form a commutative subsemigroup which is a semilattice with respect to the natural order, given by $e\leq f \Leftrightarrow ef=e.$ 

\smallskip

Let $T$ be a commutative ring. For a $T$-module  $M$ denote by  ${\rm Ann}_T(M)$   the  annihilator of $M$  in $T$.
If $M$ is a finitely generated  $T$-module,  the sets ${\rm Supp}_T (M)=\{\mathfrak{p}\in {\rm Spec}(T)\,|\,
M_\mfp\ne 0\}$ and $V({\rm Ann}_T(M))=\{\mathfrak{p}\in {\rm Spec}(T)\,|\, \mathfrak{p}\supseteq {\rm Ann}_T(M)\}$
coincide, (see e.g. \cite[p. 25-26]{HM}).

\vu

The following lemma characterizes the idempotents of ${\bf PicS} (R).$

\vm

\begin{lema} \label{idemp} Let $M$ be a f.g.p. $R$-module and $I_M={\rm Ann}_R(M)$. Then, the following statements are equivalent:

\vu

\begin{itemize}

\item[(i)] $M\ot M\cong M$.

\vu

\item[(ii)] $M\cong R/{I_M}$.

\vu

\item[(iii)] $M\cong Re$ (and $I_M=R(1-e)$), for some idempotent $e$ of $R$.
\end{itemize}
\end{lema}

\p (i)$\Rightarrow$(ii)
It easily follows from the dual basis lemma  that $M$ is a faithfully projective $\left(R/I_M\right)$-module.
Moreover,  the $R$-module isomorphism $M\ot M\cong M,$ implies $M\ot_{R/I_M} M\cong M$ as $R/I_M$-modules,
and hence ${\rm rk}_{R/I_M}(M)\le 1$.
Moreover,
${\rm Supp}_{R/I_M} (M)=V(\bar{0})={\rm Spec}(R/I_M),$
and since $M$ is a faithfully projective $R/I_M$-module, then    $[M] \in {\bf Pic} (R/I_M).$ Being an idempotent,  $[M]$ must be the identity element of   ${\bf Pic} (R/I_M), $ so that   there is an $R/I_M$-module
isomorphism $M\cong R/I_M, $  which is clearly an isomorphism of  $R$-modules.

\vd

(ii)$\Rightarrow$(iii) Since $M\cong R/{I_M}$ is f.g.p. as an $R$-module, then the exact sequence $$0\to I_M\to R\to R/{I_M}\to 0$$
 splits,  thus the ideal $I_M$ is a direct summand of  $R$ and the assertion easily follows.

\vd

(iii)$\Rightarrow$(i) It is clear.\cua

\vt

Let ${\bf I_p}(R)$  be the semilattice of the idempotents of $R$ with respect to the product.    If  the $R$-modules   $Re \cong   Rf$ are isomorphic where $e,f \in {\bf I_p}(R)$, then their annihilators in $R$ coincide, i.e. $(1-e)R= (1-f)R.$ This yields   $e=f,$ and   it follows by Lemma \ref{idemp}  that the map $e \mapsto [eR]$ is an isomorphism of ${\bf I_p}(R)$ with  the semilattice of the idempotents of ${\bf PicS} (R).$ Consequently, the  components of ${\bf PicS} (R)$ can be
indexed by the idempotents of $R,$ and $${\bf PicS} (R)  = \bigcup\limits_{e\in {\bf I_p}(R)}{\bf Pic S}_e (R)$$
gives the decomposition of ${\bf PicS} (R)$ as a semilattice of abelian groups.   Thus, if for each $e\in {\bf Ip}(R)$ we denote by   $[M_e]$  the identity element of ${\bf PicS }_e(R),$ then $[M_e][M_f]=[M_{ef}]$ and   ${\bf PicS} _e(R){\bf PicS} _f(R)\subseteq {\bf PicS}_{ef} (R),$ for all $e,f\in {\bf I_p}(R)$ (this also can be seen directly  from Lemma \ref{idemp}).

\vd

Now  we will describe the components of ${\bf PicS} (R).$ For this, note that for any $R$-module $N$ we have ${\rm Ann}_R(N)={\rm Ann}_R({\rm End}_R(N)),$
and if $N$ is projective it follows from the dual basis lemma that ${\rm Ann}_R(N)={\rm Ann}_R(N^*).$

\begin{lema} \label{dec} ${\bf PicS} _e(R)=\{[N]\in {\bf PicS} (R)\mid{\rm Ann}_R(N)=R(1-e)  \} \cong  {\bf Pic} (Re),$  for all $e\in {\bf Ip}(R).$  In particular, the identity element of ${\bf PicS} _e(R)$ is $[Re].$
\end{lema}
\p For the first equality let $[N]\in {\bf PicS} _e(R).$ Then, there are $R$-module isomorphisms ${\rm End}_R(N)\cong N\ot N^*\cong M_e\cong Re,$  where the last isomorphism follows from Lemma \ref{idemp}, and 
 hence ${\rm Ann}_R(N)={\rm Ann}_R({\rm End}_R(N))=R(1-e).$ Conversely, if ${\rm Ann}_R(N)=R(1-e)$ we have $N\ot  M_e\cong N\ot Re\cong Ne=Ne\oplus N(1-e)=N,$
as $R$-modules, so $[N]\in{\bf PicS} _e(R). $  Thus for any $[N]\in{\bf PicS} _e(R),$ its representative is a faithfully projective $Re$-module and the isomorphism ${\bf PicS} _e(R) \cong  {\bf Pic} (Re)$ is now trivial.   Finally,  notice that the image of $e$ in  $R_\mfp $ is either $0$ (if  $e \in \mfp$) or the identity of
$R_\mfp $  (if $e \notin \mfp$). Hence   $ (e R )_\mfp $ is  free of rank $0$ or $1,$ and,  consequently,   $[eR] \in {\bf PicS}(R).$ 

\vd

Summarizing, we have.
\begin{teo} \label{picde} The (disjoint) union
\begin{equation}\label{pics}{\bf PicS} (R)  = \bigcup\limits_{e\in {\bf I_p}(R)}{\bf Pic S}_e (R)\cong\bigcup\limits_{e\in {\bf I_p}(R)}{\bf Pic} (Re)
\end{equation}  gives the decomposition of   ${\bf PicS} (R)$ as a semilattice of abelian groups, whose structural homomorphisms are given by  
$\e _{e,f}:   {\bf Pic} (Re) \to {\bf Pic} (Rf) ,$ $[M] \mapsto [M \otimes Rf],$ where $e,f  \in {\bf I_p}(R),$ $ e \geq f.$ \cua
\end{teo}

We point out  the following.
\begin{lema} \label{iguald}For any $g\in G,$ we have:

\begin{itemize}

\item[(i)] ${\bf PicS}_{1_g}(R)\cong {\bf Pic}(D_g).$

\vu

\item[(ii)] Let $g\in G$ and $[M]\in {\bf PicS}(R).$  If $1_gm=m,$ for all $m\in M$ and $M_\mfp\cong (D_g)_\mfp$ as $R_\mfp$-modules, for all $\mfp\in{\rm Spec}(R),$ then   $[M]\in {\bf Pic}(D_g).$

\item [(iii)] ${\bf PicS} (D_g)  = \bigcup\limits_{e\in {\bf I_p}(R)\atop e1_g=e}{\bf Pic} (Re),$ for any $g\in G .$ 
\end{itemize}
\end{lema}

\p Item (i)  is clear from Theorem \ref{picde}.

\vu

(ii)  Notice that $M$  is a f.g.p.  $D_g$-module. Let $\mfp\in {\rm Spec}(D_g).$  Since we have a ring isomorphism $D_g\cong R/{\rm Ann}_R(D_g) ,$  we may consider $M$ as an $R/{\rm Ann}_R(D_g)$-module and make the identification $\mfp=\bar{\mfp}_1=\mfp_1/{\rm Ann}_R(D_g),$ where $\mfp_1\in{\rm Spec}(R)$ and $\mfp_1$ contains ${\rm Ann}_R(D_g).$ Thus, it follows from the assumption that there exist a $(D_g)_{\bar\mfp_1}$-module isomorphisms $M_\mfp\cong M_{{\bar \mfp}_1}\cong (D_g)_{\bar\mfp_1}\cong  (D_g)_\mfp$, which imply $[M]\in {\bf Pic}(D_g).$

\vu

(iii) Since \eqref{pics} holds for any commutative ring, we get $${\bf PicS} (D_g)  = \bigcup\limits_{e\in {\bf I_p}(D_g)}{\bf Pic} (D_ge)=\bigcup\limits_{e \in {\bf I_p}(D_g)}{\bf Pic} (Re).$$ Moreover,  $e\in {\bf I_p}(D_g)$  exactly when $e $ is an element of  ${\bf I_p}(R)$ and $e1_g=e.$  \cua  

\vu

\section{ A partial action  on  {\bf PicS}($R$) and the sequence
 $$H^1(G,\af,R)\stackrel{\varphi_1}\to {\bf Pic}(R^\af)\stackrel
{\varphi_2}\to {\bf Pic}(R )\cap {\bf Pic S}(R)^{\af^*} \stackrel{\varphi_3}\to H^2(G,\af, R)$$ }
\label{phi1phi2phi3}
\subsection{A partial action  on  PicS($R$) }

Let $\af=(D_g, \afg)_{g\in G}$ be a unital partial action of a group $G$ on    $R.$ It is known that $\afg(1_h1_{g\m})=1_g1_{gh},$
for all $g,h\in G$  (see \cite[p. 1939]{DE}). Then for any $y\in R,$ we have
\begin{equation}\label{prodp}\afg(\afh(y1_{h\m})1_{g\m})=\af_{gh}(y1_{(gh)\m})1_g, \,\,\,\text{ for all}\,\, g,h\in G. \end{equation}

\vu

In  all what follows  $\af=(D_g, \afg)_{g\in G}$ will be a fixed unital partial action of the group $G$ on  the ring  $R.$ The next result will
help us in the construction of a partial action on ${\bf PicS}(R).$

\begin{lema}\label{gaction} Let E and F be central $R$-$R$-bimodules and $g\in G.$ Suppose that $1_{g\m}x=x$ and $1_{g\m}y=y$ for all
$x\in E$ and $y\in F.$ Denote by $E_g$ the set  $E$ where  the (central) action of $R$ is given by $$r\bullet x_g=
\af_{g\m}(r1_g)x=x_g\bullet r,\,\, r\in R,\,\, x_g\in E_g.$$  Then

\vu

\begin{itemize}
\item[(i)]  $E_g$ is an $R$-module and $(E_g)_\mfp=(E_\mfp)_g$ as $R$-modules,  where the action of $R$ on $(E_\mfp)_g$ is
$r\bullet \frac{x}{s}=\frac{\af_{g\m}(r1_g)x}{s},$ for any $x\in E,\,\mfp \in {\rm Spec}(R), \,s\in R\setminus \mfp.$

\vd

\item[(ii)]  ${\rm Hom}_R(E,F)= {\rm Hom}_R(E_{g},F_{g })$ as sets. In particular, we have ${\rm Iso}_R(E,F)=
{\rm Iso}_R(  E_{g },F_{g } )$ and  ${\rm End}_R(E) = {\rm End}_R(E_{g } ).$

\vd

\item[(iii)] If $E$ is a f.g.p. $R$-module, so too is $E_g.$

\vd

\item[(iv)] There is an R-module isomorphism $(E\ot F)_g\cong E_g\ot F_g.$

\vd

\item[(v)] If $\rm{rk} ($E$)\le 1,$ then $\rm{rk} ({\it E_g})\le 1.$

\vd

\item[(vi)] For any $[M]\in {\bf Pic}(D_{g\m}),$ $[M_g]\in{\bf Pic}(D_{g}). $
\end{itemize}
\end{lema}

\p  Item (i) is clear.

\vd

(ii) Obviously  ${\rm Hom}_R(E,F)\subseteq {\rm Hom}_R( E_{g},F_{g }  ).$ Let $f\in {\rm Hom}_R(  E_{g },F_{g }),\, r\in R$ and $x\in E_{g } .$
Then $f(rx)=f(r1_{g^{-1}} x)=f(\af_{g\m}(r'1_{g})x)=f(r'\bullet x)=r'\bullet f(x)=rf(x),$ where $r'\in R$ is such that $\af_{g\m}(r'1_{g})=r1_{g\m}.$

\vd

(iii)  For any $f\in E^*,$ the map
$$\afg\circ f\colon E_g \ni x\mapsto \afg \circ f(x)=\afg(f(x)1_{g\m})\in R$$  is an element of $(E_g)^*.$ Indeed,
$$\afg \circ f(r\bullet x)=\afg(f(r\bullet x)1_{g\m})=\afg(\af_{g\m}(r1_{g})f(x))=r(\afg \circ f)(x).$$

Suppose that $E$ is  a f.g.p. $R$-module. Then, there are $f_i\in E^*$ and $x_i\in E$ such that
$x=\sum\limits_{i}f_i(x)x_i=\sum\limits_{i}f_i(x)1_{g\m}x_i=\sum\limits_{i}(\afg\circ f_i(x))\bullet x_i,$
for any $x\in E$, which implies that $E_g$ is a f.g.p. $R$-module, with dual basis $\{\af_g\circ f_i, x_i\}.$

\vd

(iv)  The map $E_g\times F_g\ni (x,y)\mt (x\ot y)_g\in (E\ot F)_g$ is $R$-balanced, therefore it induces a well defined $R$-module map
\begin{equation}\label{iotag}\iota_g\colon E_g\ot_R F_g\ni x\ot y \mt (x\ot y)_g \in(E\ot F)_g\end{equation}
which is clearly bijective.

\vd

(v) Take  $\mfp \in {\rm Spec}(R).$ We have two cases to consider.

\emph{Case 1:} $E_\mfp=0.$ In this case we have $(E_g)_\mfp\stackrel{{\rm (i)}}{=}(E_\mfp)_g=0.$

\emph{Case 2:} $E_\mfp\cong R_\mfp$ as $R_\mfp$-modules. Since $1_{g\m}x=x,$ for all $x\in E$ we have
$1_{g\m} r_\mfp= r_\mfp$ for all $r_\mfp\in R_\mfp,$ which implies $1_{g\m}\notin\mfp,$ because the
image of 
 $1_{g\m}$ in  $R_\mfp$ is either $0$ or the identity of  $R_\mfp $, and  thus
$E_\mfp\cong R_\mfp =  (D_{g\m})_\mfp\,\,\, \text{as}\,\,\,  R_\mfp\text{-modules.}$

Finally, using (i) we get
$$(E_g)_\mfp=(E_\mfp)_g\cong((D_{g\m})_\mfp)_g=((D_{g\m})_g)_\mfp\cong  (D_{g})_\mfp ,$$
where the latter isomorphism is given by $\alpha _g.$ This ensures that  $\rm{rk}({\it E_g})\le 1.$

\vd

vi) In the proof of  item (iii) we saw that if $\{f_i, m_i\}$ is a dual basis for $M,$ then $\{\af_g\circ f_i, m_i\}$ is a dual basis
for the $D_g$-module $M_g.$ On the other hand, by the same reason as in the proof of item (ii),  we have that
$D_{g} \cong D_{g\m}\cong {\rm End}_{D_{g\m}}(M)={\rm End}_{D_{g}}(M_g)$ as rings. Since $M$ is faithful, so too is $M_g,$ and the
ring isomorphism $D_g\cong {\rm End}_{D_{g}}(M_g)$ implies  $[M_g]\in {\bf Pic}(D_{g}). $ \cua

\begin{lema}\label{igual} For any $g\in G$ set $$X_g=\{[1_gE]\,|\, [E]\in \,{\bf PicS}(R)\}=[D_g]{\bf PicS}(R).$$

Then, $X_g$ is an ideal of ${\bf PicS} (R)$ and

\begin{itemize}

\item[(i)] $X_g=\{[E]\in {\bf PicS} (R)\,|\, E=1_gE\},$

\vu

\item[(ii)] For any $[E]\in X_{g\m}$ we have  $[E_g]\in X_g.$
\end{itemize}
\end{lema}
\p  (i) It is clear that $X_g\supseteq\{[E]\in \,${\bf PicS}(R)$\,|\, E=1_gE\}.$ On the other hand, given
$[E]\in X_g $ there exists $[F]\in {\bf PicS} (R)$ and an $R$-module isomorphism $\varphi_g\colon 1_gF\to E.$
This leads to $E=\varphi_g(1_gF)=1_g\varphi_g(1_gF)=1_gE.$

\vd

(ii) Notice that $1_g\bullet x_g=1_{g\m}x_g=x_g$ for any $x_g \in E_g,$ so $1_g\bullet E_g=E_g.$ By
Lemma \ref{gaction} $E_g$ is a f.g.p. $R$-module and ${\rm rk}(E_g) \le 1,$ hence $[E_g] \in {\bf PicS}(R).$
Thus, using item (i) we conclude that $[E_g]\in X_g.$ \cua

\vd

\begin{teo} \label{psem} The family $\af^*=( X_g, \afg^*)_{g\in G},$ where $\afg^*\colon X_{g\m}\ni [E]\mapsto [E_g]\in X_{g}$
defines a partial action of G on ${\bf PicS} (R).$
\end{teo}
\p By Lemmas \ref{gaction} and \ref{igual} the map  $\afg^*$ is a well defined semigroup homomorphism, for all $g\in G.$
Clearly $X_1={\bf PicS} (R)$ and  $\af_1^*={\rm id}_{{\bf PicS} (R)}.$ We need to show that $\af_{gh}^*$ is an extension
of $\afg^*\circ \afh^*.$ If $[E]\in X_{h\m}$ is such that $[E_h]\in X_{g\m},$ then $E=1_{h\m}E$ and
$$E=E_h=1_{g\m}\bullet E_h=\af_{h\m}(1_{g\m}1_h)E_h=1_{(gh)\m}1_{h\m}E=1_{(gh)\m}E.$$
Thus $E =1_{(gh)\m}E,$  which shows that ${\rm dom}(\afg^*\circ \afh^*)\subseteq {\rm dom}\,\af^*_{gh},$ thanks to item (i)
of Lemma \ref{igual}. Furthermore, we have $\afg^*\circ \afh^*([E])=[(E_h)_g],$\, $\af^*_{gh}([E])=[E_{gh}],$ and
$(E_h)_g=E_{gh}$ as sets. Now, for any $r\in R,\, x=(x_h)_g\in (E_h)_g$ we get

\begin{align*}
r\bullet(x_h)_g &=\af_{h\m}(\af_{g\m}(r1_g)1_h)x
\stackrel{(\ref{prodp})} {=}\af_{(gh)\m}(r1_{gh})1_{h\m}x
=\af_{(gh)\m}(r1_{gh})x
=r\bullet 1_{gh}x,
\end{align*}
and $ (E_h)_g\cong E_{gh}$ as $R$-modules. In particular,  $\afg^*$ has an inverse $\af_{g\m}^*$ , so that
each $\afg^*$ is an isomorphism. \cua

\begin{remark}\label{ggaction} It follows from Theorem \ref{psem} and {\rm (\ref{prodp})} that  there is an $R$-module isomorphism
$$( D_{(gh)\m}\ot  P)_{gh}\ot D_g\cong ( D_{g\m}\ot (D_{h\m}\ot P)_h)_g,$$ for any $R$-module $P$ and $g,h\in G.$
\end{remark}

\vu

The subset of invariants of ${\bf PicS}(R)$ (see equation (\ref{inva})) is given by
\begin{equation}\label{fixx}{\bf PicS}(R)^{\af^*}=\{[E]\in {\bf PicS}(R)\,|\, (D_{g\m}\ot E)_g\cong D_g\ot E, \,\,\text{for all} \,g\in G\}.
\end{equation}

\vu

\begin{prop}  ${\bf PicS}(R)^{\af^*}$ has an element 0 and is a commutative  inverse submonoid of ${\bf PicS}(R).$
\end{prop}
\p Evidently $0\in {\bf PicS}(R)^{\af^*}.$ Moreover, for any $[E],[N]\in {\bf PicS}(R)^{\af^*},$ we have
$\af^*_g([D_{g\m}\ot (E\ot N)])=\af^*_g([D_{g\m}\ot E])\af^*_g([D_{g\m}\ot  N)])=[D_g\ot E][D_g\ot N]=[D_g\ot (E\ot N)]$
and $[E][N]\in {\bf PicS}(R)^{\af^*}. $

\vu

Given any element $[E]$ of  ${\bf PicS}(R)^{\af^*},$ we need to show that  $[E^*]$ is also in ${\bf PicS}(R)^{\af^*}.$
If $[E]\in X_{g\m},$ for some $g\in G,$ then $[E^*]=[E^*][E][E^*]\in X_{g\m}.$
Since $[E^*][E][E^*]=[E^*]$ and  $[E ][E^*][E ]=[E ],$ then $[(E^*)_g][E_g][(E^*)_g]=[(E^*)_g]$ and $[E_g][E^*_g][E_g]=[E_g],$
thanks to  (iv) of Lemma \ref{gaction} . Thus,
\begin{equation}\label{star}[(E^*)_g]=[E_g]^*=[(E_g)^*],\end{equation}
where the last equality follows from  Proposition \ref{picinv}. Therefore for any $[E]\in  {\bf PicS}(R)^{\af^*}$ we get
\begin{align*}\{\af^*_g([D_{g\m}\ot E^*])\}^* &=[(D_{g\m}\ot E^*)_g]^*=[((D_{g\m}\ot  E)^*)_g]^* \stackrel{\eqref{star}}=[(D_{g\m}\ot E)_g]\\
&=[D_g\ot E]=[D_g\ot E]^{**}=[D_g\ot E^*]^*,\end{align*}
hence $\af^*_g([D_{g\m}\ot  E^*])=[D_g\ot E^*],$ and $[E^*]\in {\bf PicS}(R)^{\af^*}.$ Finally, since $\afg^*$ is
a ring isomorphism, we have $\afg^*([D_{g\m}])=[D_g],\, g\in G,$ or equivalently $[R]\in {\bf PicS}(R)^{\af^*}.$ \cua

\begin{remark} \label{unxg}
Recall that by definition $X_g={\bf PicS}(R)[D_g],$ and thus by Theorem \ref{picde} we have  that $X_g=\bigcup\limits_{e\in {\bf I_p}(R)}{\bf Pic} (Re) [D_g]=\bigcup\limits_{e\in {\bf I_p}(R)}{\bf Pic} (D_ge),$  where the last equality holds because the map ${\bf Pic} (Re)\ni [M]\mapsto [M\ot D_ge ]\in {\bf Pic} (D_ge) $ is a group epimorphism. Consequetly,  $$X_g=\bigcup\limits_{e\in {\bf I_p}(R)}{\bf Pic} (D_ge)=\bigcup\limits_{e\in {\bf I_p}(R)}{\bf Pic} (Re1_g)=\bigcup\limits_{e\in {\bf I_p}(D_g)\atop e1_g=e}{\bf Pic} (Re)={\bf PicS}(D_g),$$ thanks to {\rm (iii)} of Lemma \ref{iguald}. Finally,  Remark~\ref{unit} implies  that $\U(X_g)={\bf Pic}(D_g).$
\end{remark}

\vt

In all what follows  $G$ will stand for  a finite group and  $R\supseteq R^\af$ will be  an $\af$-partial Galois extension.  In particular, it follows from \cite[(ii) Theorem 4.1]{DFP} that  $R$ and any $D_g, g\in G,$ are f.g.p $R^\af$-modules.
We also recall two well known results that will be used  several times in the sequel. The proof of the first can be easily obtained by localization.
For the second one we refer to the literature.

\begin{lema}\label{proje} Let  $M,N$ be R-modules such that $M$ and  $M\ot N$ are f.g.p. R-modules, and $M_\mfp\ne 0$ for all $\mfp \in {\rm Spec}(R).$
Then, N is also a  f.g.p.  R-module.
\end{lema}

\vu

\begin{lema}\label{cancel}\cite[Chapter I, Lemma 3.2 (b)]{DI} Let  $M$ be a faithfully projective  R-module. Then, there is a $R$-$R$-bimodule
isomorphism $M^*\ot_{{\rm End}_R(M)} M\cong R.$ Consequently if $N,N'$ are $R$-modules such that  $M\ot N\cong M\ot N'$ as $R$-modules,
then $N\cong N'$ as $R$-modules.\cua
\end{lema}

\vu

\subsection{The  sequence $H^1(G,\af,R)\stackrel{\varphi_1}\to{\bf Pic}(R^\af)\stackrel{\varphi_2}
\to {\bf PicS}(R)^{\af^*}\cap {\bf Pic}(R)\stackrel{\varphi_3}\to H^2(G,\af,R)$}

For any   $R\star_\af G$-module $M,$ as in \cite[page 82]{DFP} we denote $$M^G=\{m\in M\,|\, (1_g\de_g)m=1_gm, \,\,\text{for all}\,\, g\in G\}.$$
It  can be seen that  $R$ is an  $R\star_\af G$-module via $(r_g\de_g)\rhd r=r_g\afg(r1_{g\m}),$
for each $g\in G$ and $r\in R,r_g\in D_g,$ and this action induces an $R\star_\af G$-module structure on $R\ot_{R^\af} M^G.$

\vu

We have the following.

\begin{teo}\label{homp} There is a group homomorphism $\varphi_1\colon H^1(G,\af, R) \to {\bf Pic}(R^\af).$
\end{teo}
\p  Take $f\in Z^1(G,\af, R)$ and define $\te_f\in {\rm End}_{R^\af}(R\star_\af G)$ by $\te_f(r_g\de_g)=r_gf(g)\de_g$ for all
$r_g\in D_g, g\in  G.$ Then, $\te_f$ is an $R^\af$-algebra homomorphism because
\begin{align*}
\theta _f(r_g\de_g)\theta _f (r_h\de_h)&=(r_gf(g)\de_g)(r_hf(h)\de_h)=r_gf(g)\af_g(r_hf(h)1_{g\m})\de_{gh}\\
&=r_g\af_g(r_h1_{g\m})f(g)\af_g(f(h)1_{g\m})\de_{gh}
=r_g\af_g(r_h1_{g\m})f(gh)\de_{gh}\\
&=\te_f((r_g\de_g)(r_h\de_h)),
\end{align*} for all $r_g\in D_g, r_h\in D_h, g,h\in G.$ Hence, we may define an $R\star_\af G$-module $R_f$ by $R_f=R$ as sets and
 \begin{equation*}\label{rf}(r_g\de_g)\cdot r=\te_f(r_g\de_g) \rhd r,\,\,\,\text{for any}\,\, r\in R, \, r_g\in D_g,\, g\in G.\end{equation*}
 In particular $R_f=R$ as $R$-modules,  as  $f$  is  normalized in view of Remark~\ref{normalized}.  By (iii) of \cite[Theorem 4.1]{DFP} there is an $R$-module isomorphism $R\ot_{R^\af}R_f^G \cong R_f$. 
Since $R$ is a f.g.p. $R^\af$-module  we conclude that  $R_f^G $ is a  f.g.p. $R^\af$-module by Lemma \ref{proje}.
 Finally, via localization we see from the last isomorphism  that ${\rm rk}_{R^\af}(R_f^G)=1,$ so $[R_f^G]\in {\rm Pic}(R^\af).$ Define
\begin{equation*}\label{1homo}\varphi_1\colon H^1(G,\af, R) \ni {\rm cls}(f) \to [R_f^G]\in  {\rm Pic}(R^\af).\end{equation*}
We will check that $\varphi_1$ is a well defined group isomorphism. If $f\in B^1(G,\af, R),$ there exists $a\in \U(R)$ such that
$f(g)=\afg(a1_{g\m})a\m,$ for all $g\in G.$ In this case one has $ r\in R_f^G$ if and only if,  $ar\in R^\af.$
Thus, multiplication by $a$ gives an $R^{\af}$-module isomorphism $R_f^G\cong R^\af,$ which yields
$[R_f^G]= [R^\af] \in {\rm Pic}(R^\af).$ Therefore, to prove that $\f_1$ is well defined, it is enough to show that
$\varphi_1$ preserves products. For any $f,g \in Z^1(G,\af, R)$ there is a chain of $R^\af$-module isomorphisms

$$ R\ot _{R^{\af}} (R_f^G\ot_{R^\af} R_g^G)\cong (R\ot_{R^\af}  R_f^G)\ot (R\ot_{R^\af}
{R_g}^G)\cong R_f \ot R_g \cong R \ot R   \cong R_{fg} \cong R\ot _{R^{\af}} R_{fg}^G,$$

and recalling that $R$ is a f.g.p. $R^\alpha$-module we have $R_f^G\ot _ { R^\af  }R_g^G \cong R_{fg}^G$  as $R^\af$-modules, by
Lemma \ref{cancel}.    \cua

\begin{prop} There is a group homomorphism $ {\bf Pic}(R^\af)\stackrel{\varphi_2}{\to }{\bf PicS}(R)^{\af^*}\cap {\bf Pic}(R).$
\end{prop}
\p For any $[E]\in  {\bf Pic}(R^\af)$ set $\varphi_2([E])=[R\ot_{R ^\af} E].$ Clearly  $\varphi_2$ is a well defined
group homomorphism from ${\bf Pic}(R^\af)$ to ${\bf Pic}(R)$. We shall check that
${\rm im}\,\varphi_2\subseteq {\bf PicS}(R)^{\af^*}.$

\vu

There are $R$-module isomorphisms
\begin{equation}\label{priso}D_g\ot(R\ot_{R ^\af} E)\cong D_g\ot_{R^\af}E\,\,\text{ and}\,\,
(D_{g\m  } \ot (R\ot_{R ^\af} E))_{g}\cong(D_{g\m} \ot_{R^\af} E  )_{g} .\end{equation}
   Furthermore,  the map determined by
\begin{equation}\label{seiso} D_g\ot_{R^\af}E \ni d\ot_{R^\af} x\mapsto \af_{g\m}(d) \ot_{R^\af} x \in (D_{g\m} \ot_{R^\af} E )_g,\end{equation}
 is also an  $R$-module isomorphism.  Then, combining \eqref{priso} and \eqref{seiso} we obtain an $R$-module isomorphism
 $D_{g}\ot (R\ot_{R^\af} E)\cong (D_{g\m}\ot (R\ot_{R^\af} E))_{g},$ for all $g\in G,$  and hence $\varphi_2([E])\in {\bf PicS}(R)^{\af^*}.$ \cua

\vd

Now, we proceed with the construction of  $\varphi_3.$ First, for any $R$-module $M$ we identify $M\ot D_g$ with $MD_g,$ via
the $R$-module isomorphism $M\ot  D_g\ni x\ot  d\mt xd \in MD_g$,  for any $g\in G.$  Now,  let $[E]\in  {\bf PicS}(R)^{\af^*}\cap {\bf Pic}(R).$
Then, by \eqref{fixx} there is a family of $R$-module isomorphisms
\begin{equation}\label{psiii}\{\psi_g\colon ED_g\to (ED_{g\m})_g\}_{g\in G} \;\;  \text{with}\;\;  \psi_g(rx)=\af_{g\m}(r1_g)\psi_g(x), \end{equation}
where $r\in R, x\in ED_g, g\in G.$
Thus the maps $\psi\m_g\colon (ED_{g\m})_g\to ED_g,\,\, g\in G,$ satisfy \begin{equation}\label{psii}\psi\m_g(rx)=\psi
\m_g ( \afg(r1_{g\m})\bullet x)=\afg(r1_{g\m})\psi\m_g(x),\,\,\text{ for all}\,\,r\in R, \,\,x \in ED_{g\m}.
\end{equation}
We shall prove that  $\psi_{(gh)\m}\psi\m_{h\m}\psi\m_{g\m}\colon ED_gD_{gh}\to  ED_gD_{gh}$ is well defined and is an element
of $\U({\rm End}_{D_gD_{gh}}(ED_gD_{gh})).$ From (\ref{psii}) we have $\psi\m_{g\m}(ED_gD_{gh})
\subseteq ED_{g\m}D_{h}\subseteq {\rm dom}\,\psi\m_{h\m}\cap\, {\rm dom}\,\psi_{g\m},$ and also
$\psi\m_{h\m}(ED_{g\m}D_{h})\subseteq ED_{h\m g\m}D_{h\m}\subseteq {\rm dom}\,\psi_{h\m g\m}\cap \, {\rm dom}\,\psi_{h\m}.$
This yields  that the map $\psi_{(gh)\m}\psi\m_{h\m}\psi\m_{g\m}$ is well defined and   $\psi_{g\m}\psi_{h\m}\psi\m_{(gh)\m}$ is its  inverse.

Now we check that $\psi_{(gh)\m}\psi\m_{h\m}\psi\m_{g\m}$ is $D_gD_{gh}$-linear. Take $d\in D_gD_{gh}$ and $x\in ED_gD_{gh}.$
Using (\ref{psii}) and (\ref{psiii}) we get the following
\begin{align*}
&\psi_{(gh)\m}\psi\m_{h\m}\psi\m_{g\m}(dx) \stackrel{d\in D_gD_{gh}}{=}  \psi_{(gh)\m}(\af_{h\m} ( \af_{g\m}(d)))\psi\m_{h\m}\psi\m_{g\m}(x) =\\
&\psi_{(gh)\m}(\af_{h\m g\m}(d))\psi\m_{h\m}\psi\m_{g\m}(x)= d\psi_{(gh)\m}\psi\m_{h\m}\psi\m_{g\m}(x).
\end{align*}
Since $[E]\in  {\bf Pic}(R) $ then $  [ED_gD_{gh}]\in {\bf Pic}(D_gD_{gh})$, and thus ${\rm End}_{D_gD_{gh}}(ED_gD_{gh})\cong D_gD_{gh}.$
Moreover,  $\psi_{(gh)\m}\psi\m_{h\m}\psi\m_{g\m}$ is an invertible element of
${\rm End}_{D_gD_{gh}}(ED_gD_{gh}),$ and hence  there exists $\om_{g,h}\in \U(D_gD_{gh})$ such that $\psi_{(gh)\m}\psi\m_{h\m}\psi\m_{g\m}(x)=\om_{g,h}x,$
for all $x\in ED_gD_{gh}.$

Summarizing, for an element $[E]\in {\bf PicS}(R)^{\af^*}\cap {\bf Pic}(R) $ we have found a map
$\om_{[E]}=\om\colon G\times G\ni (g,h)\mapsto \om_{g,h}\in \U(D_gD_{gh})\subseteq R.$
We shall see that $\om \in Z^2(G,\af,R).$ Take $g,h,l\in G$ and $x\in ED_gD_{gh}D_{ghl}.$ Then
\begin{align*}
\om_{g,hl}\afg(\om_{h,l}1_{g\m})x\stackrel{\eqref{psiii}}=&\om_{g,hl}\psi_{g\m}(\om_{h,l}\psi\m_{g\m}(x))\\
=&\psi_{(ghl)\m}\psi\m_{(hl)\m}\psi\m_{g\m}\psi_{g\m}(\om_{h,l}\psi\m_{g\m}(x))\\
=&\psi_{(ghl)\m}\psi\m_{(hl)\m}(\om_{h,l}\psi\m_{g\m}(x))\\
=&\psi_{(ghl)\m}\psi\m_{(hl)\m}(\psi_{(hl)\m}\psi\m_{l\m}\psi\m_{h\m}\psi\m_{g\m}(x))\\
=&\psi_{(ghl)\m}\psi\m_{l\m}\psi\m_{h\m}\psi\m_{g\m}(x)\\
=&\psi_{(ghl)\m}\psi\m_{l\m}\psi\m_{(gh)\m}\psi_{(gh)\m}\psi\m_{h\m}\psi\m_{g\m}(x)\\
=& \om_{gh,l}\om_{g,h}x.
\end{align*}
Notice that $ [E D_gD_{gh}D_{ghl}] \in {\bf Pic} ( D_gD_{gh}D_{ghl} )$ and,  in particular,  $E D_gD_{gh}D_{ghl}$ is a
faithful $D_gD_{gh}D_{ghl}$-module.  Since $\om_{g,hl}\afg(\om_{h,l}1_{g\m}),  \om_{gh,l}\om_{g,h} \in D_gD_{gh}D_{ghl},$ we obtain that
$\om_{g,hl}\afg(\om_{h,l}1_{g\m})= \om_{gh,l}\om_{g,h}$ as desired.

\vd

\begin{claim}\label{IsoChoice}  cls$(\om)$ does not depend on the choice of the isomorphisms.\end{claim}

Let $\{\lb_g\,|\, g\in G\}$ be another choice of $R$-isomorphisms $ ED_g\to (ED_{g\m})_g.$ Then $\lb_g(rx)=\af_{g\m}(r1_g)\lb_g(x),$
$\lb\m_g(rx)=\afg(r1_{g\m})\lb\m_g(x)$, for all $g\in G,\,r\in R$ and $x$ belonging to the correspondent domain. Let 
$\tilde\om\colon G\times G\to R$ also be defined by
$$\tilde\om(g,h)x=\lb_{(gh)\m}\lb\m_{h\m}\lb\m_{g\m}(x),\,\,\,\text{ for any}\,\,\,x\in ED_gD_{gh}.$$
We shall prove that ${\rm cls}(\om)={\rm cls}(\tilde\om)$ in  $H^2(G,\af,R).$
 Since $\lb_g\psi\m_g\colon( ED_{g\m})_g\to (ED_{g\m})_{g}$ is $D_{g}$-linear and  $[(ED_{g\m})_g]\in {\bf Pic}(D_g),$ there exists $u_g\in \U(D_g)$
 such that $\lb_g\psi\m_g$ is the multiplication by $u_g.$ Then, the map $u\colon G\ni g\to u_g\in R$ belongs to $C^1(G,\af,R)$,
 and for any $x\in ED_gD_{gh}D_{ghl}$ we have
\begin{align*}
\om\m_{g,h}\tilde\om_{g,h}x&=\psi_{g\m}\psi_{h\m}\psi\m_{(gh)\m}(\tilde\om_{g,h}x)\\
&=(\psi_{g\m}\lb\m_{g\m})\lb_{g\m}(\psi_{h\m}\lb\m_{h\m})\lb\m_{g\m}(\lb_{g\m}\lb_{h\m}\psi\m_{(gh)\m})(\tilde\om_{g,h}x)\\
&=u\m_{g\m}(\lb_{g\m}u\m_{h\m}\lb\m_{g\m})(\lb_{g\m}\lb_{h\m}\psi\m_{(gh)\m})(\!\!\!\underbrace{\tilde\om_{g,h}}_{\in\, D_gD_{gh}}\!\!\!\!x)\\
&=u\m_{g\m}(\lb_{g\m}u\m_{h\m}\lb\m_{g\m})\tilde\om_{g,h}(\lb_{g\m}\lb_{h\m}\psi\m_{(gh)\m})(x)\\
&=u\m_{g\m}(\lb_{g\m}u\m_{h\m}\lb\m_{g\m})(\lb_{(gh)\m}\lb\m_{h\m}\lb\m_{g\m}\lb_{g\m}\lb_{h\m}\psi\m_{(gh)\m})(x)\\
&=u\m_{g\m}(\lb_{g\m}u\m_{h\m}\lb\m_{g\m})(\lb_{(gh)\m}\psi\m_{(gh)\m})(x)\\
&=u\m_{g\m}(\lb_{g\m}u\m_{h\m}\lb\m_{g\m})(u_{(gh)\m}x)\\
&=u\m_{g\m}\lb_{g\m}[u\m_{h\m}\af_{g\m}(u_{(gh)\m}1_{g})\lb\m_{g\m}(x)]\end{align*}\begin{align*}
&=u\m_{g\m}\afg(u\m_{h\m}1_{g\m})u_{(gh)\m}x\\
&=v_{g}\afg(v_{h}1_{g\m})v\m_{(gh) }x,
\end{align*}
where  $v_{g}=u\m_{g\m}.$  Since the map $v\colon G\ni g\to v_g\in R$ belongs to $C^1(G,\af,R),$ this shows that
${\rm cls}(\om)={\rm cls}(\tilde\om)$ in  $H^2(G,\af,R)$ as desired.

\vt

Let $\varphi_3\colon {\bf PicS}(R)^{\af^*}\cap {\bf Pic}(R)\ni [E] \mapsto {\rm cls(\om)}\in H^2(G,\af, R).$

\begin{claim}  $\varphi_3$ does not depend on the choice of the representative of $[E],$ for any $[E]\in{\bf PicS}(R)^{\af^*}\cap  {\bf Pic}(R).$
\end{claim}

Let $[E]=[F]\in {\bf PicS}(R)^{\af^*}\cap {\bf Pic}(R),$ and $\{\psi_g\,|\, g\in G\},\,\{\lb_g\,|\, g\in G\}$ be families
of $R$-isomorphisms $ED_g\to (ED_{g\m})_g,$ $FD_g\to (FD_{g\m})_g$ inducing  the $(2,\af)$-cocycles $\om,\tilde{\om}$ respectively.
Let also $\Omega\colon E\to F$ be an $R$-module isomorphism. Then $\Omega(ED_g)=\Omega(E)D_g=FD_g$ and one obtains the
$R$-module isomorphisms $\Omega|_{ED_g}\colon ED_g\to FD_g$  and  $\Omega\mid _{ED_{g\m}}\colon (ED_{g\m})_g \to (FD_{g\m})_g$
(thanks to (ii) of Lemma \ref{gaction}), for any $g\in G.$  Thus, the family $\{\Om\psi_g\Om\m\mid_{FD_g}\colon FD_g\to (FD_{g\m})_g\,|\, g\in G\}$
induces a $(2,\af)$-cocycle which is cohomologous to $\tilde{\om},$ in view of Claim~\ref{IsoChoice}, and  we may suppose that $\Om\psi_g\Om\m=\lb_g$
on $FD_g,g\in G.$ Hence, given $x\in FD_gD_{gh}$ we have
\begin{align*}\tilde{\om}_{g,h}x=\lb_{(gh)\m}\lb\m_{h\m}\lb\m_{g\m}(x)\equiv \Om\psi_{(gh)\m}\psi_{h\m}\psi_{g\m}\Om\m(x)=\Om(\om_{g,h}\Om\m(x))=\om_{g,h}x,
\end{align*} which implies $\tilde{\om}_{g,h} = {\om}_{g,h},$ because $ FD_gD_{gh}$, being an element of  ${\bf Pic}(D_gD_{gh})$, is a faithful  $D_gD_{gh}$-module.
We conclude  that in general ${\rm cls}(\om)={\rm cls}(\tilde\om)$ and $\varphi_3$ is well defined.

\begin{claim} $\varphi_3$ is a group homomorphism.\end{claim}
Let $[E],[F]\in {\bf PicS}(R)^{\af^*}\cap{\bf Pic}(R)$ with $\varphi_3([E])={\rm cls}(\om)$ and $\varphi_3([F])={\rm cls}(\om').$
Consider  families of $R$-module isomorphisms $\{\phi_g\colon ED_g\to (ED_{g\m})_g\}_{g\in G}$ and $\{\lb_g\colon FD_g\to (FD_{g\m})_g\}_{g\in G}$
defining ${\rm cls}(\om)$ and ${\rm cls}(\om')$ respectively.

Notice that $ED_g\ot FD_g=(E\ot F)D_g,$ for all $g\in G,$ and by  (ii) and (iv) of Lemma \ref{gaction},
$(ED_{g\m})_g\ot (FD_{g\m})_g\cong ((E\ot F)D_{g\m})_g, \, g\in G,$ via the map $\iota_g$ defined in \eqref{iotag}. Then,
$$\iota_g\circ (\phi_g\ot \lb_g) \colon (E\ot F)D_g\to ((E\ot F)D_{g\m})_g,\, g\in G $$ is an $R$-module isomorphism which induces an element $u\in Z^2(G,\af, R),$ where
\begin{align*}
u_{g,h} x\ot  y &=[\iota_{(gh)\m}\!\circ\! (\phi_{(gh)\m}\ot \lb_{(gh)\m} ) ][\iota_{h\m}\circ\!(\phi_{h\m}\ot \lb_{h\m}) ]\m[\iota_{g\m}\!\circ (\phi_{g\m}\ot \lb_{g\m})  ]\m (x\ot  y)\\
&=(\phi_{(gh)\m}\phi \m _{h\m}\phi \m _{g\m}\ot \lb_{(gh)\m} \lb \m _{h\m} \lb \m _{g\m} )(x\ot  y)\\
&=(\om_{g,h}\ot \om'_{g,h})(x\ot  y)=(\om_{g,h}\om'_{g,h})(x\ot  y),
\end{align*}
for all $x\in ED_gD_{gh}, y\in FD_gD_{gh}, \, g,h\in G.$ Then, $u=\om\om'$ and $\varphi_3$ is a group homomorphism.

\vu

\section{The sequence $H^2(G,\af,R)\stackrel{\f_4}\to B(R/R^\af)\stackrel{\f_5}\to H^1 (G,\af^*,{\bf PicS}(R))$ }\label{phi4phi5}
We start this section by giving  some preliminary  results that help us to construct the homomorphism $\f_4$.

\vu

First of all we recall from \cite{DES} that
the partial crossed product $R\star_{\af,\om}G$ for the unital twisted partial action  $(\af,\om)$ of $G$ on $R$ is the direct sum $\bigoplus_{g\in G}D_g\de_g$,
in which the $\de_g'$s are symbols, with the multiplication defined by the rule:
$$(r_g\de_g) � (r'_h\de_h) = r_g\af_g(r'_h1_{g\m})\om_{g,h}\de_{gh},$$ for all $g,h\in G$, $r_g\in D_g$ and $r'_h\in D_h$. If, in particular, the twisting $\om$
is trivial, then we recover the partial skew group ring $R\star_\af G$ as given in Subsection 2.3.

\begin{prop} \label{iso}If $\om, \tilde{\om}\in Z^2(G,\af, R)$ are cohomologous, there is an isomorphism of $R^\af$-algebras and
$R$-modules $R\star_{\af,\om} G\cong R\star_{\af,\tilde{\om}} G.$
\end{prop}
\p There exists $u\in C^1(G,\af, R),$ $u\colon G\ni g\mapsto u_g\in \U(D_g)\subseteq R$ such that $\om_{g,h}=\tilde{\om}_{g,h} u_g\af_g(u_h1_{g\m})u\m_{gh} ,$
for all $g,h\in G.$  Take  $a_g\in D_g,$  set $\f(a_g\de_g)=a_gu_g\de_g$  and extend $\f$ to $\f\colon R\star_{\af,\om} G\to R\star_{\af,\tilde{\om}} G$
by $R$-linearity. Clearly $\f$ is bijective with inverse $a_g\de_g\mapsto a_gu\m_g\de_g, $  then we only need to prove that $\f$ preserves products.
\begin{align*}
&\f((a_g\de_g)(b_h\de_h))=\f(a_g\afg(b_h1_{g\m})\omgh\de_{gh})=
a_g\afg(b_h1_{g\m})\omgh u_{gh}\de_{gh}=\\
&a_g\afg(b_h1_{g\m})\tilde{\om}_{g,h}u_g\afg(u_h1_{g\m}) \de_{gh}=(a_gu_g\de_{g})(b_hu_h\de_{h}) =\f(a_g\de_g)\f(b_h\de_h).
\end{align*} \cua

We also give the next.

\begin{prop}\label{rpiso}  If $\f\colon R\star_{\af,\om} G\to R\star_{\af} G$ is an isomorphism of $R^\af$-algebras and $R$-modules,
then $\om$ is cohomologous to  $1_{\af}=\{1_g1_{gh}\}_{g,h\in G}.$
\end{prop}
\p To avoid confusion, we write $R\star_{\af,\om} G=\bigoplus\limits_{g\in G}D_g\de_g,$
$R\star_{\af} G=\bigoplus\limits_{g\in G}D_g\de'_g$ and identify $R=R\de_1=R\de'_1.$ For $r\in R$ we have $\f(r\de_1)=r\de'_1\in R$ and
$$(1_g\de_g)r(\om\m_{g\m,g}\de_{g\m})
=\afg(r1_{g\m})\af_g(\om\m_{g\m,g})\om_{g,g\m}\de_{1}
\stackrel{(\ref{afom})}{=}\afg(r1_{g\m})\de_1.$$
Hence, $\f(1_g\de_g)r\f(\om\m_{g\m,g}\de_{g\m})=\f(\afg(r1_{g\m})\de_1)=\afg(r1_{g\m})\de'_1=(1_g\de'_g)r(1_{g\m}\de'_{g\m}).$ From this we obtain $$(1_{g\m}\de'_{g\m})\f(1_g\de_g)r\f(\om\m_{g\m,g}\de_{g\m})=(1_{g\m}\de'_{g\m})(1_g\de'_g)r(1_{g\m}\de'_{g\m})=r(1_{g\m}\de'_{g\m}),$$
and, multiplying by the right both  sides of the equality by $\f(1_g\de_g)$, we get
\begin{align*}
r[(1_{g\m}\de'_{g\m})\f(1_g\de_g)]&=[(1_{g\m}\de'_{g\m})\f(1_g\de_g)r\f(\om\m_{g\m,g}\de_{g\m})]
\f(1_g\de_g)\\
&=(1_{g\m}\de'_{g\m})\f(1_g\de_g)r\f[(\om\m_{g\m,g}\de_{g\m})(1_g\de_g)]\\
&=(1_{g\m}\de'_{g\m})\f(1_g\de_g)r\f(1_{g\m}\de_1)\\
&=(1_{g\m}\de'_{g\m})\f(1_g\de_g)\f(1_{g\m}\de_1)r\\
&=(1_{g\m}\de'_{g\m})\f((1_g\de_g)(1_{g\m}\de_1))r\\
&=(1_{g\m}\de'_{g\m})\f(1_g\de_g)r,
\end{align*}
so $(1_{g\m}\de'_{g\m})\f(1_g\de_g)\in C_{R\star_{\af } G}(R).$

\vu

Since $R$ is commutative and $R\supseteq R^\af$ is an $\af$-partial Galois extension, \cite[Lemma 2.1(vi) and Proposition 3.2]{PS} imply
 that $R\star_{\af,\tilde\om} G$ is $R^\af$-Azumaya  and $C_{R\star_{\af,\tilde\om} G}(R)= R$ for arbitrary $\tilde{\om} ,$ in particular this is
 true for $R\star_{\af } G.$ Thus, $$(1_{g\m}\delta'_{g\m})\f(1_g\de_g)=r_g,\,\, \text{for some}\,\,r_g\in R,$$ and, multiplying from the left both the sides of
 the last equality by $1_{g}\delta'_{g}$ we obtain $\f(1_g\de_g)=u_g\delta'_g,$ where $u_g=\afg(r_g1_{g\m})\in D_g.$ Therefore, $\f(1_g\de_g)= u_g\delta'_g, \, g\in G.$

\vu

On the other hand there exists $W=\sum\limits_{ h\in G}   a_h \de_h $ such that $1_g\delta'_g=\f(W)=\sum\limits_{h\in G}a_hu_h {{\de }'_h},$  then
$1_g\delta'_g=a_gu_g\delta'_g$ and $u_g\in \U(D_g).$ Set $u\colon G\ni g\mt u_g\in \U(D_g)\subseteq R,$ then $u\in C^1(G,\af, R)$ and $$\om_{g,h}u_{gh}\delta'_{gh}=\f(\om_{g,h}\de_{gh})=\f(1_g\de_g)\f(1_h\de_h)=(u_g\delta'_g)(u_h\delta'_h)=u_g\afg(u_h1_{g\m})\delta'_{gh}.$$
From this we conclude that  $\om_{g,h}u_{gh}= u_g\afg(u_h1_{g\m})$, hence $\om$ is cohomologous to $1_{\af}.$\cua

\vu

\begin{prop} \label{pb} Let $\om\in Z^2(G,\af, R).$ Then $R\star_{\af,\om} G$ is an  Azumaya  $R^\af$-algebra  and $[R\star_{\af,\om} G]\in B(R/R^\af).$
\end{prop}
\p As mentioned above, the facts that $R$ is commutative and $R\supseteq R^\af$ is an $\af$-partial Galois extension  imply that $R\star_{\af,\om} G$ is
$R^\af$-Azumaya  and $C_{R\star_{\af,\om} G}(R)= R.$ The latter implies that $R$ is a maximal commutative $R^\af$-subalgebra of $R\star_{\af,\om} G.$
On the other hand, by \cite[Theorem 4.2]{DFP} the extension $R\supseteq R^\af$ is separable, and
 finally \cite[Theorem 5.6]{AG} tells us that $R\star_{\af,\om} G$ is split by $R,$ which means that $[R\star_{\af,\om} G]\in B(R/R^\af).$\cua

\vd

It follows from  Propositions \ref{iso} and \ref{pb} that there is a well defined function
\begin{equation}\label{fun} \varphi_4\colon H^2(G,\af, R)\ni {\rm cls}(\om)\mapsto [R\star_{\af,\om} G]\in B(R/R^\af). \end{equation}
\begin{subsection}{ $\varphi_4$ is a homomorphism}

In this subsection we follow the ideas of \cite[chapter IV]{DI}  to prove that  the map $\varphi_4$ defined in (\ref{fun}) is  a group  homomorphism.
For this we need a series of lemmas.

\begin{lema} \label{iso1}Let   $\om \in   Z^2(    G, \af, R).$ Then, there exists an $R^\af$-algebra isomorphism\linebreak
$(R\star_{\af,\om}G)^{op}\cong R\star_{\af,\om\m} G.$
\end{lema}
\p  Define $\phi\colon (R\star_{\af,\om}G)^{op}\to R\star_{\af,\om\m} G$  by $\phi(r_g\delta_g)=\af_{g\m}(r_g\om_{g,g\m})\de_{g\m},$
for all $g\in G,\, r_g\in D_g.$ Note that $\phi$ is an $R^\af$-module isomorphism with inverse
$R\star_{\af,\om\m} G\ni r_{g\m}\de_{g\m} \mt \af_g(r_{g\m}\om\m_{g\m,g}) \de _g\in (R\star_{\af,\om}G)^{op} $. For $g,h\in G,r_g\in D_g$
and $t_h\in D_h$ we have
\begin{align*}
\phi[(r_g\de_g)\circ(t_h\de_h)]&=\phi(t_h\afh(r_g 1_{h\m})\om_{h,g}\de_{hg})\\
&=\af_{g\m h\m}(t_h\afh(r_g 1_{h\m})\om_{h,g}\om_{hg,g\m h\m})\de_{g\m h\m}\\
&=\af_{g\m}[\af_{h\m}(t_h\afh(r_g 1_{h\m})\om_{h,g}\om_{hg,g\m h\m})]\de_{g\m h\m}\\
&\stackrel{(\rm {v})}{=}\af_{g\m}[\af_{h\m}(t_h\afh(r_g 1_{h\m})\afh(\om_{g,g\m h\m})\om_{h, h\m})]\de_{g\m h\m}\\
&\stackrel{(\ref{afom})}{=}\af_{g\m}(\af_{h\m}(t_h)r_g \om_{g,g\m h\m}\om_{h\m, h})\de_{g\m h\m}.
\end{align*}
On the other hand
\begin{align*}
\phi(r_g\de_g)\phi(t_h\de_h)&=(\af_{g\m}(r_g\om_{g,g\m})\de_{g\m})(\af_{h\m}(t_h\om_{h,h\m})\de_{h\m})\\
&=(\af_{g\m}(r_g)\om_{g\m,g})(\af_{g\m}(\af_{h\m}(t_h)\om_{h\m,h})\om\m_{g\m,h\m}\de_{g\m h\m})\\
&=\af_{g\m}(r_g\af_{h\m}(t_h)\om_{h\m,h})\om_{g\m,g}\om\m_{g\m,h\m}\de_{g\m h\m}\\
&\stackrel{(\rm {v})}{=}\af_{g\m}(r_g\af_{h\m}(t_h)\om_{h\m,h})\af_{g\m}(\om_{g,g\m h\m})\de_{g\m h\m}\\
&=\af_{g\m}(r_g\af_{h\m}(t_h)\om_{h\m,h}\om_{g,g\m h\m})\de_{g\m h\m},
\end{align*}
and we conclude that $\phi$ is multiplicative. \cua

\vd

From now on, in order to simplify notation we will denote $R^e=R\otimes_{R^\af}R$.

\begin{lema}\label{idem} There is a family of orthogonal idempotents $e_g\in R^e$ with $g\in G,$ satisfying the following properties:

\begin{equation}\label{i} (1\otimes_{R^\af} \afg(x1_{g\m}))e_g=(x\otimes_{R^\af} 1)e_g,\,\,\text {for all} \,\,x\in R.\end{equation}

\begin{equation} \label{ii} \sum\limits_{g\in G}e_g=1_{R^e}.\end{equation}
\end{lema}
\p By (iv) of \cite[Theorem 4.1]{DFP} the map $\psi \colon R^e\to \prod_{g\in G}D_g,$ given by
$\psi(x\ot_{R^\af} y)=(x\afg(y1_{g\m}))_{g\in G},$ is an isomorphism of $R$-algebras.
Take $v_g=(x_h)_{h\in G} \in \prod_{h\in G}D_h , g\in G,$ where $x_h=\de_{h,g}1_g.$ Then, the set $\{e_g=\psi\m(v_{g\m}) |\, g\in G\}$
is a family of orthogonal idempotents in $R^e.$ Using the isomorphism $\psi$ we  check \eqref{i} and \eqref{ii}.\cua

\begin{lema}\label{iso2} For any $\om \in Z^2(G,\af, R)$  there is an $R$-module isomorphism  $R\star_{\af,\om} G\cong R^e.$
\end{lema}
\p Let  $\{e_g\,|\, g\in G\}$ be the family of pairwise orthogonal idempotents constructed in  Lemma \ref{idem}, and set
$\eta \colon R\star_{\af,\om} G\to R^e$ defined by $\eta(\sum_{g\in G}r_g\de_g)=\sum_{g\in G}(r_g\ot_{R^\af} 1)e_{g\m}.$
Clearly $\eta$ is $R$-linear and  we only need to check that $\eta$ is an isomorphism.

If $\sum_{g\in G}(r_g\ot_{R^\af} 1)e_{g\m}=0,$ then  $(r_h\ot_{R^\af} 1)e_{h\m}=0,$ for any $h\in G.$ Applying the isomorphism $\psi$
from the proof of Proposition \ref{idem} we conclude that $r_h=0$ and $\eta$ is injective.

Now we prove the surjectivity. By applying  $\psi$ we get  $(1\ot_{R^\af} r1_{g\m})e_{g\m}=(1\ot_{R^\af}r )e_{g\m},$ $r\in R, g\in G.$
Then for any $r,s\in R,$ we obtain
\begin{align*}r\ot_{R^\af}s&\stackrel{\eqref{ii}}=\sum_{g\in G}(r\ot_{R^\af}s)e_{g\m}=\sum_{g\in G}(r\ot_{R^\af}s1_{g\m})e_{g\m}\\
&\stackrel{\eqref{i}}=\sum_{g\in G}(r\af_g(s1_{g\m})\ot_{R^\af}1)e_{g\m}=\eta \left(\sum_{g\in G}r\af_g(s1_{g\m})\de_g\right).
\end{align*} \cua
\medskip

By  Lemma \ref{iso2} the map $\bar\eta\colon {\rm End}_{R^\af}(R\star_{\af,\om}G)\ni f\mapsto \eta f\eta\m\in {\rm End}_{R^\af}(R^e),$
is an $R^\af$-algebra isomorphism.

\vu

Now we prove that the tensor product of partial Galois extensions is also a partial Galois extension. More precisely we have the following.
\begin{prop} \label{tensorgal}Let $G$ and $H$ be finite groups and $\af=(D_g, \afg)_{g\in G},\0=(I_h, \0_h)_{h\in H}$ (unital)
partial actions of $G$ and $H$ on  commutative rings $R_1$ and $R_2,$ respectively. Assume that $R_1^\af=R_2^\0=k$ and suppose that the ring extensions
$R_1\supseteq k,$ $R_2\supseteq k$ are $\af$- and $\0$-partial Galois, respectively. Then $R_1\otimes_k R_2\supseteq k \otimes_kk=k$ is
an $(\af\otimes_k\0)$-partial Galois extension.
\end{prop}
\p In this proof  unadorned $\otimes$ means  $\otimes_k.$  Note that $\af\otimes \0=( D_g\otimes I_h, \afg\otimes \0_h)_{(g,h)\in G\times H},$
is partial action of $G\times H$ on the ring $R_1\ot R_2.$ Consider $x=\sum_{i}u_i\otimes v_i\in k\ot k= R_1^\af\otimes R_2^ \0,$ then for
all $(g,h)\in G\times H$ we have $\afg\otimes\0_h(x (1_{g\m}\otimes 1_{h\m}))=\sum_{i}u_i1_{g}\otimes v_i1_{h}=x(1_g\otimes 1_h),$ and
$(R_1\otimes R_2)^{\af\otimes \0}\supseteq R_1^\af\otimes R_2^ \0.$

\vu

Conversely, take  $c_1\in R_1$ and $c_2\in R_2$ such that  ${\rm tr}_{R_1/R_1^\af}(c_1)=1_{k}={\rm tr}_{R_2/R_2^\0}(c_2).$ Now let
$x\in (R_1\otimes R_2)^{\af\otimes\0}$ and write $x(c_1\otimes c_2)=\sum\limits_{i=1}^{m}t_i\otimes s_i.$ Then,
\[
\begin{array}{ccl}
x&=&x(1_{R_1}\otimes 1_{R_2})=\sum_{(g,h)\in G\times H}x(\afg(1_{g\m}c_1)\otimes \0_h(1_{h\m}c_2))\\
&=& \sum_{(g,h)\in G\times H}\!x(1_g\!\otimes \!1_h)(\afg\!\otimes\!\0_h)(1_{g\m}c_1\!\otimes\!1_{h\m}c_2)\\
&=&\sum_{(g,h)\in G\times H}\!(\afg \otimes\0_h)(x(c_1\!\otimes\! c_2)(1_{g\m}\!\otimes\!1_{h\m}))\\
&=&\sum_{i=1}^{m}\sum_{(g,h)\in G\times H}(\afg\otimes\!\0_h)(t_i1_{g\m}\otimes s_i1_{h\m})\\
&=&\sum_{i=1}^{m}\sum_{(g,h)\in G\times H}\afg(t_i1_{g\m})\otimes\0_h( s_i1_{h\m})\\
& =&\sum_{i=1}^{m}{\rm tr}_{R_1/R_1^\af}(t_i)\otimes{\rm tr}_{R_2/R_2^\0}(s_i)\in R_1^\af\otimes R_2^ \0 .
\end{array}
\]
Thus  $(R_1\otimes R_2)^{\af\otimes \0}=R_1^\af\otimes R_2^ \0=k.$

\vu

 Finally, let $m,n\in \N$ and $\{x_i, y_i\,|\, 1\le i\le m\},$ $\{u_i, z_i\,|\, 1\le i\le n\}$ be the $\af$- and $\0$-partial Galois
 systems for the extensions $R_1\supseteq k$ and $R_2\supseteq k,$ respectively. Then for all $(g,h)\in G\times H$  we see that
\begin{align*}
&\sum_{i,j}(x_i\otimes u_j)(\afg\otimes\0_h)((y_i\otimes z_j)(1_{g\m}\otimes1_{h\m}))=\sum_{i,j}(x_i\otimes u_j)(\afg(y_i1_{g\m})\otimes\0_h(z_j1_{h\m}))\\
&=\sum_{i}x_i\afg(y_i1_{g\m})\otimes \sum_ju_j\0_h(z_j1_{h\m})=\de_{1,g}\otimes \de_{1,h}=\de_{(1,1), (g,h)}. \end{align*}

We conclude that the set  $\{x_i\otimes u_j, y_i\otimes z_j\,|\, 1\le i\le m, 1\le j\le n\}$ is an $\af\otimes\0$-partial Galois coordinate system for the
extension $R_1\otimes R_2\supseteq k \otimes k.$ \cua

\vd

With the same hypothesis and notations given in Proposition \ref{tensorgal}, we have.
\begin{prop} \label{isotens}There is an isomorphism of $k$-algebras
$$(R_1\star_{\af,\om} G)\ot _k (R_2\star_{\0,\tilde{\om}} H)\cong (R_1\ot_{k}R_2)\star_{\af \ot \0,\om\ot \tilde{\om}}(G\times H).$$
\end{prop}
\p Here again  unadorned $\otimes$ means  $\otimes_k.$ We denote
$$S=(R_1\ot R_2)\star_{\af \ot\0,\om\ot\tilde{\om}}(G\times H)=\bigoplus\limits_{(g,h)\in G\times H}(D_g\otimes I_h)\ep_{(g,h)},$$
$S_1=R_1\star_{\af,\om} G=\bigoplus\limits_{g\in G}D_g\delta_g$ and  $S_2=R_2\star_{\0,\tilde{\om}} H={\bigoplus\limits_{h \in H }I_h \delta'_h}.$
For $(g,h)\in G\times H,$ the map
$$D_g\de_g\times I_h\de'_h\ni (a_g\de_g,b_h\de'_h)\mapsto (a_g\ot b_h)\ep_{(g,h)}\in S$$ extended by $k$ linearity to $S_1\times S_2$ is clearly a bilinear  $k$-balanced map.
 Hence, it induces a bijective $k$-linear map  $\xi\colon S_1\ot S_2\to S$ such that $a_g\de_g\ot b_h\de'_h\mapsto (a_g\ot b_h)\ep_{(g,h)}.$
 The fact that  $\xi$ preserves products is straightforward.
\cua
\begin{prop} \label{2co} Given $\om \in Z^2(G,\af, R) ,$ we have $\om\ot_{R^\af} \om\m\in B^2(G\times G,\af\ot _{R^\af}\af, R^e).$ Thus if
$\tilde\om\in Z^2(G,\af, R), $ then $\om\ot_{R^\af}\tilde\om$ is cohomologous to $\om\tilde\om \ot_{R^\af} 1_{\af}.$
\end{prop}
\p By Proposition \ref{tensorgal}, $R^e$ is an $\af\ot_{R^\af} \af$-partial Galois extension of $R^\af$ with the partial action of
$G\times G.$ Then using the isomorphisms appeared in Proposition \ref{isotens}, Lemmas \ref{iso1}, \ref{iso2}, \cite[Theorem 2.1 (c)]{AG} and iv) of
\cite[Theorem 4.1]{DFP} we obtain  a chain of $R^\af$-algebra isomorphisms
\begin{align*}
&(R^e)\star_{\af\ot_{R^\af}\af, \om\ot_{R^\af}\om\m}(G\times G)\stackrel{}{\to}(R\star_{\af,\om} G)\ot_{R^\af}(R\star_{\af,\om\m} G) \stackrel{{\rm id}\ot_{R^\af}\phi\m}{\to}\\
&(R\star_{\af,\om} G)\ot_{R^\af}(R\star_{\af,\om} G)^{\rm{op}}\stackrel{\Gamma}{\to}\, {\rm End}_{R^\af}(R\star_{\af,\om} G)\stackrel{\bar\eta}{\to}\,{\rm End}_{R^\af}(R^e)\stackrel{j\m}{\to}
\\
&\,(R^e)\star_{\af\ot\af}(G\times G),
\end{align*}
where $\Gamma$ is given by $\Gamma(x\ot y):z\mapsto xzy$, for all $x,y,z\in R\star_{\af,\om} G$, and $j$ is given by \eqref{jota}.
By Proposition \ref{rpiso} one only needs to show that the above composition  restricted to  $R^e$ is the identity. We have
\begin{align*}
(r_1\ot_{R^\af} t_1)\de_{(1,1)}&\mapsto r_1\de_1\ot_{R^\af} t_1\de_1\mapsto y=r_1\de_1\ot_{R^\af} t_1\de_{1}\mapsto \Gamma(y)\mapsto j\m(\eta \Gamma(y)\eta\m).
\end{align*}

 For a fixed $a\in G$ we compute the image of $\eta \Gamma(y)\eta\m$ on $(r_a\ot_{R^\af} 1)e_{a\m}\in R^e,\, r_a\in D_a.$ We see that $\eta \Gamma(y)\eta\m((r_a\ot_{R^\af} 1)e_{a\m})$
\begin{align*}
=&\eta \Gamma(y)(r_a \de_a)=\eta((r_1\de_1)(r_a\delta_a)(t_1\de_{1}))=\eta(r_1r_a\af_{a}(t_11_{a\m })\de_{a })\\
=&(r_1r_a\af_{a}(t_11_{a\m})\ot_{R^\af} 1 ) e_{a\m}
=(r_1r_a\ot_{R^\af} 1)(\af_{a}(t_11_{a\m})\ot_{R^\af} 1)e_{ a\m }\\
\stackrel{\eqref {i}}{=}&(r_1r_a\ot_{R^\af} 1)(1\ot_{R^\af} t_11_{a\m})e_{a\m }=
(r_1r_a\ot_{R^\af} t_1)(1\ot_{R^\af} \af_{a\m}(1_a))e_{a\m }\\
\stackrel{\eqref{i}}{=}&(r_1\ot_{R^\af} t_1)(r_a\ot_{R^\af} 1)e_{a\m }=j((r_1\ot_{R^\af} t_1)\de_{(1,1)})[(r_a\ot_{R^\af} 1)e_{a\m }].
\end{align*}

Then, $\eta \Gamma(y)\eta\m(x)=j((r_1\ot_{R^\af} t_1)\de_{(1,1)})x,$ for all $x\in R^e,$  and we conclude that
$$j\m(\eta \Gamma(y)\eta\m) =(r_1\ot_{R^\af} t_1)\de_{(1,1)}.$$
Hence, the composition is $R^e$-linear and    $\om\ot_{R^\af} \om\m\in B^2(G\times G,\af\ot_{R^\af} \af, R^e).$

\vu

Finally, for any  $\tilde\om\in {  Z^2(G ,  \af, R)},$ we have $(\om\ot_{R^\af}\tilde\om)(\tilde\om\ot_{R^\af}\tilde\om\m)=
\om\tilde\om\ot_{R^\af} 1_\af,$ and  this yields that $\om\ot_{R^\af}\tilde\om$ is cohomologous to $\om\tilde\om\otimes_{R^\af}1_\af.$ \cua

\begin{teo} \label{homo}Let R be an $\af$-partial Galois extension of $R^\af.$ Then the map \linebreak
$ \varphi_4\colon H^2(G,\af, R)\ni {\rm cls}(\om)\mapsto [R\star_{\af,\om} G]\in B(R/R^\af) $ is a group homomorphism.
\end{teo}
\p In this proof  unadorned $\otimes$ will mean   $\otimes_{R^\af}.$  Let ${\rm cls}(\om), {\rm cls}(\tilde \om)\in H^2(G,\af, R). $
By  Propositions \ref{isotens}, \ref{iso} and \ref{2co}  we have
\begin{align*}
[R\star_{\af,\om} G][R\star_{\af,\tilde\om} G]&=[(R\ot R)\star_{\af\ot \af,\om\ot\tilde\om} (G\times G)]=[(R\ot R)\star_{\af\ot \af,\om\tilde\om\ot 1_{\af}} (G\times G)]\\
&=[R\star_{\af,\om\tilde\om} G][R\star_{\af,1_{\af}} G]=[R\star_{\af,\om\tilde\om} G][{\rm End}_{R^\af}(R)]=[R\star_{\af,\om\tilde\om} G].
\end{align*}  which gives $[R\star_{\af,\om} G][R\star_{\af,\tilde\om} G]=[R\star_{\af,\om\tilde\om} G]$ and the assertion follows.\cua\end{subsection}

\vu

\subsection{The construction of $B(R/R^\af)\stackrel{\f_5}\to H^1 (G,\af^*,{\bf PicS}(R)) $}
We remind that  unadorned $\ot$ stands for $\ot _R.$

\vu

Let $\af^*=(\af^*_g, X_g)_{g\in G}$ be  the partial action of $G$ on ${\bf PicS}(R)$ constructed  in Theorem \ref{psem}.  Since 
$\U( {\bf PicS}(R))={\bf Pic}(R) $ we have that
$B^1(G,\af^*, {\bf PicS}(R))$ is the group $$\{ f\in C^1(G,\af ^*,{\bf PicS}(R)) \mid f(g) = \af^*_g([P] [D_{g^{-1}}])  [P^{*}], \, \mbox{for some}\,\, [P]\in  {\bf Pic}(R)  \}$$
and $	Z^1(G,\af^*, {\bf PicS}(R))$ is given by
 \begin{align*}
	 \{f\in C^1 (G,\af ^*,{\bf PicS}(R))\mid f(gh)[D_g]=f(g) \, \af^* _g ( f(h) [D_{g^{-1}}] ),\,\, \forall g,h\in G \}.\end{align*}

\begin{remark}\label{trivi}Let $f\in C^1 (G,\af,{\bf PicS}(R)), \, g\in G$ and $\mfp \in {\rm Spec}(R).$ We shall make a little abuse of
notation by writing $f(g)_\mfp$ for a representative of the class $f(g)$ localized at $\mfp.$ 
\end{remark}

\medskip
We proceed with the construction of $\f_5.$ Take $[A]\in B(R/R^\af).$ Then by \cite[Theorem 5.7]{AG} there is an Azumaya $R^\af$-algebra
equivalent to $A$  containing $R$ as a maximal commutative subalgebra. Hence, we assume that $A$ contains $R$ as a maximal commutative subalgebra.
By \cite[Theorem 4.2]{DFP} $R\supseteq R^\af$ is separable, moreover \cite[Theorem 5.6]{AG} tells us that $A$ is a faithfully projective
$R$-module and  there is a $R$-algebra isomorphism $R\ot_{R^\af} A^{\rm{op}}\cong {\rm End}_R(A).$

On the other hand,
$D_g\ot A$ is a  faithfully projective $D_g$-module, thus  by Proposition \ref{ht} there is an $R$-algebra isomorphism
${\rm End}_{D_g}(D_g\ot _R A)\cong D_g\ot {\rm End}_R(A),$ for any $g\in G.$ Consequently, we have an $R$-algebra isomorphism
${\rm End}_{D_g}(D_g\ot  A)\cong D_g\ot_{R^\af} A^{\rm{op}}.$

 Therefore, by \cite[Proposition I.3.3]{DI} the functor \begin{equation}\label{cequiv}\_\_\ot_{D_g} (D_g\ot A)
 \colon _{D_g}{\rm Mod}\to _{(D_g\ot_{R^\af} A^{\rm{op}})}\!\!{\rm Mod},\end{equation}  determines a category equivalence.

It is clear that  $D_g\ot A$ is a  left $R\ot_{R^\af}A^{\rm{op}}$-module via $(r\ot_{R^\af} a)(d'\ot a')=rd'\ot a'a.$ Moreover,
let $(D_{g\m}\ot A)^g = D_{g\m}\ot A$ (as sets) endowed with  a left $R\ot_{R^\af}A^{op}$-module structure via
\begin{equation} \label{ngaction} (r\ot_{R^\af} a)\bullet (d'\ot a')=\af_{g\m}(r1_g)d'\ot a'a,\end{equation} for any
{ $g\in G, r\in R, a\in A,  d'\in D_{g\m}.$}
Restricting, we obtain left $D_g\ot_{R^\af}A^{op}$-module structures on $D_g\ot A$ and $(D_{g\m}\ot A)^g,$ respectively.
Moreover  $D_{g\m}\ot A$ is  also a right $R\ot_{R^\af} A^{op}$-module via
\begin{equation}\label{ngaction1}(d\ot a)(r\ot_{R^\af} a')= dr\ot a'a.\end{equation}

Furthermore, we denote by $_{h}(D_{g\m}\ot A)_I,$ the $R$-$R$-bimodule $D_{g\m}\ot A,$ where the actions of $R$ are induced
by (\ref{ngaction}) and  (\ref{ngaction1}), and $h\in\{g,1_G\}.$

\vu

It follows from  (\ref{ngaction}) that  $(D_{g\m}\ot A)^g$ is an object in ${}_{D_g\ot_{R^\af} A^{\rm{op}}}{\rm Mod}$   and 
by \eqref{cequiv} there is a  $D_g$-module $M^{(g)}$ such that
\begin{equation}\label{unico}(D_{g\m}\ot A)^g\cong M^{(g)}\ot_{D_g} (D_g\ot A) \cong M^{(g)} \ot A,\,
\,\text{as}\,\, (D_g\ot_{R^\af} A^{\rm{op}})\text{-modules},\end{equation}
 where $M^{(g)}$ is considered as an  $R$-module via the map $r\mt r1_g,\, r\in R.$ Our aim is to show
 that $[M^{(g)}]\in  {\bf PicS}(R),$ for all $g\in G.$
 
\vu

From (\ref{unico}) we see that
\begin{equation}\label{n2}(D_{g\m}\ot A)^g\cong M^{(g)} \ot A\,\,\, \,\, \text{as} \,\, \,\,R\ot_{R^\af} A^{\rm{op}}\,\text{-modules.}\end{equation}

As $R$-modules we have $(D_{g\m}\ot A)^g=(D_{g\m}\ot A)_g.$ Since $D_{g\m}$ is a f.g.p. $R$-module, we have that $D_{g\m}\ot A$ is a
f.g.p. $R$-module too, and by  (iii) of Lemma \ref{gaction} we conclude that $(D_{g\m}\ot A)^g$ is also a  f.g.p. $R$-module.
Since $A_\mfp\ne 0,$ for all $\mfp \in {\rm Spec}(R),$ Lemma \ref{proje} and \eqref{n2} imply that  $M^{(g)}$ is a f.g.p. $R$-module.

Now we prove that ${\rm rk}_\mfp(M^{(g)}_\mfp)\le 1,$ for all $\mfp \in {\rm Spec}(R).$  Since $A_\mfp \cong R_\mfp^{n_\mfp},$
for some $n_\mfp\ge 1$,  then  $(D_{g\m}\ot  A)_\mfp \cong ({D_{g\m}})_\mfp^{n_\mfp}.$ By (i) of Lemma \ref{gaction} we get
$((D_{g\m}\ot A)_g)_\mfp \cong (D_{g})_\mfp^{n_\mfp},$  which using \eqref{n2} implies
$$(M^{(g)}_\mfp)\ot_{R_\mfp}A_\mfp \cong  (D_{g})_\mfp^{n_\mfp} \cong  (D_{g})_\mfp\ot_{R_\mfp}{R_\mfp}^{n_\mfp}\cong (D_{g})_\mfp  \ot_{R_\mfp}A_{\mfp}.$$
We conclude that
\begin{equation}\label{rigual}{\rm rk}_\mfp(M^{(g)}_\mfp)={\rm rk}_\mfp((D_{g})_\mfp)\le1, \,\,\, \,\text{for any} \,\,\,\mfp \in {\rm Spec}(R),\end{equation}
and by Proposition \ref{carpic} we have $[M^{(g)}]\in {\bf PicS}(R).$  In addition, since $M^{(g)}$ is a (unital) $D_g$-module,
we also have that  $[M^{(g)}]\in X_g=[D_g]{\bf PicS}(R).$
\smallskip

Set $f_A\colon G\ni g\mapsto [M^{(g)}]\in {\bf PicS}(R).$ Notice that $f=f_A$ is well defined by Lemma \ref{proje}  and
$M^{(1)}=R$ satisfies (\ref{unico}). 
We shall check that $f\in Z^1(G,\af^*, {\bf PicS}(R)).$
Using (\ref{n2}), Remark \ref{ggaction} and (iv) of Lemma \ref{gaction} we obtain  $R$-module isomorphisms
\begin{align*}
(D_g\ot M^{(gh)})\ot A&\cong  D_g  \ot (M^{(gh)}\ot  A) \cong D_g\ot (D_{(gh)\m}\ot A)_{gh} \\
&\cong [D_{g\m}\ot (D_{h\m}\ot A)_{h}]_g
\cong [ D_{g\m}\ot (M^{(h)} \ot A)]_g\\
&\cong
  (D_{g\m}\ot M^{(h)}  )_g\ot( D_{g\m}\ot A )_g\cong (D_{g\m}\ot M^{(h)}  )_g\ot  M^{(g)}\ot A.
\end{align*}
Finally, by Lemma \ref{cancel} we get
\begin{equation}\label{1cocy}M^{(gh)}\ot D_g\cong (M^{(h)} \ot D_{g\m})_g\ot M^{(g)},\end{equation} which gives
$f(gh)[D_g]=f(g)\af^*_g(f(h)[D_{g\m}]).$ Taking $h=g\m$ in (\ref{1cocy}) we obtain that $D_g\cong (M^{(g\m)} \ot D_{g\m})_g\ot M^{(g)},$
as $R$-modules. Thus, $[M^{(g)}]\in \U(X_g)$ and we conclude that $f\in Z^1(G,\af^*,\! {\bf PicS}(R)).$
\medskip

 We define $\varphi_5 \colon B(R/R^\af)\ni [A]\mt {\rm cls}(f_A)\in H^1(G,\af^*, {\bf PicS}(R)).$


 \begin{claim}  $\varphi_5$  is well defined.\end{claim} Suppose $[A]=[B]\in B(R/R^\af),$ where $A$ and $B$ contain $R$ as a maximal
 commutative subalgebra (see \cite[Theorem 5.7]{AG}). There are faithfully projective $R^\af$-modules $P,Q$ such that
$A \otimes_{R^\af} {\rm End}_{R^\af}(P)\cong B \otimes_{R^\af} {\rm End}_{R^\af}(Q), $ as $R^\af$-algebras.
 It is proved in \cite[page 127]{DI} that this leads to the existence of a f.g.p. $R$-module $N$ with ${\rm rk }(N)=1$ satisfying

\begin{equation}\label{clasic}(A\ot_{R^\af} P^*)\ot N\cong B\ot_{R^\af} Q^* \,\,\, \text{as}\,\,\, R \text{-modules.} \end{equation}

We know that there are $D_g$-modules  $M^{(g)}, W^{(g)}$ such that
\begin{equation}\label{mequiv1}
 (D_{g\m}\ot A)^g\cong M^{(g)} \ot_ R A\,\,\,\, \text{as}\,\, \,\, R\ot_{R^\af} A^{\rm{op}} \text{-modules}
\end{equation}
 and
\begin{equation}\label{mequiv2}
(D_{g\m}\ot B)^g\cong W^{(g)} \ot_ R B\,\,  \,\, \text{as} \,\,\,\, R\ot_{R^\af} B^{\rm{op}}\text{-modules},
\end{equation}
for each $g\in G.$  Let $f_A,f_B\colon G\to {\bf PicS}(R)$ be defined by $$f_A(g)=[M^{(g)}], \,\,\,\,\,\,f_B(g)=[W^{(g)}],\, \,\,\,\,\,\,g\in G.$$
We must show that ${\rm cls}(f_A)={\rm cls}(f_B)$ in $H^1(G,\af^*, {\bf PicS}(R)  ).$

\vu

 Since for any $R^\af$-module $P$ one has  $(D_{g\m}\ot A \ot_{R^\af} P^*)_g\cong (D_{g\m } \ot A)_g \ot_{R^\af} P^*$
as $R$-modules, there are $R$-module isomorphisms
\begin{align*}
[D_{g\m}\ot (B\ot_{R^\af} Q^* )]_g&\stackrel{(\ref{clasic})}{\cong} [D_{g\m}\ot (A\ot_{R^\af} P^*)\ot N]_g
\cong [D_{g\m}\ot N\ot (A\ot_{R^\af} P^*) ]_g \\
&\cong (N\ot D_{g\m})_g  \ot ( D_{g\m} \ot A)_g \ot_{R^\af} P^*
\stackrel{(\ref{mequiv1})}{\cong} (N\ot D_{g\m})_g\ot M^g\ot  A \ot_{R^\af} P^*.
\end{align*} On the other hand,
\begin{equation}\label{isoB}
[D_{g\m}\ot (B\ot_{R^\af} Q^* ) ]_g \cong(D_{g\m}\ot B)_g\ot_{R^\af} Q^*
 \stackrel{(\ref{mequiv2})}{\cong} W^{(g)} \ot B\ot_{R^\af} Q^*.
\end{equation}

But $[N]\in {\bf Pic}(R),$ and,  as $R$-modules, one has
\begin{align*}
&[(D_{g\m}\ot N )_g\ot N^*\ot M^{(g)}]  \ot (N\ot  A\ot_{R^\af} P^*)\cong [(D_{g\m}\ot N )_g\ot  M^{(g)}\ot A)]\ot_{R^\af} P^*\\
&\cong [D_{g\m}\ot (B\ot_{R^\af} Q^* ) ]_g
 \stackrel{(\ref{isoB})} \cong  W^{(g)} \ot B\ot_{R^\af} Q^*
 \stackrel{(\ref{clasic})}{\cong} W^{(g)} \ot(A\ot_{R^\af} P^*)\ot N\\
&\cong W^{(g)} \ot (N\ot A\ot_{R^\af} P^*).
\end{align*}
Since $A$ is a faithfully projective $R$-module and $P^*$ is a  faithfully projective $R^\af$-module, then
$A\ot_{R^\af} P^*$ is a faithfully projective $R=R\ot_{R^\af}R^\af$-module. Therefore, $N\ot A\ot_{R^\af} P^*$
is a also a faithfully projective $R$-module, and by Lemma \ref{cancel}
$$(D_{g\m}\ot N )_g\ot N^*\ot M^{(g)}\cong W^{(g)}$$ as $R$-modules, which is equivalent to say that
$f_B(g)=\afg^*([N][D_{g\m}])[N]^{\m}f_A(g).$ This shows that $\varphi_5$ is well defined.

\medskip
\begin{teo} $\varphi_5$  is a group homomorphism.\end{teo}
\p Let $[A_1],[A_2]\in B(R/R^\af)$ and suppose that $R$ is  a maximal commutative subalgebra of $A_i, \, i=1,2.$
Let $B=A_1\ot_{R^\af}A_2.$ By \cite[Theorem 4.2]{DFP} the extension $R\supseteq R^\af$ is separable, so let
$e\in R^e$ be a separability idempotent for $R.$   Then ${\rm End}_{R^\af}(Be)$ is an Azumaya $R^\af$-algebra with
$C_{{\rm End}_{R^\af}(Be)} (B)={\rm End}_{B}(Be)\cong (eBe)^{\rm{op}},$ via the $R^\af$-algebra map $f\mt ef(e).$
It follows  from \cite[Theorem II.4.3]{DI}  that  $(eBe)^{\rm{op}}$ is an Azumaya $R^\af$-algebra and $B\ot_{R^\af}(eBe)^{\rm{op}}\cong {\rm End}_{R^\af}(Be)$
as $R^\af$-algebras. We conclude that $[B]=[eBe] $ in $B(R/R^\af).$ Following the procedure given in   \cite[page 128]{DI}
we also check that $R$ is a maximal commutative $R^\af$-subalgebra of $eBe.$ Thus, $[A_1][A_2]=[eBe]$ and there are $D_g$-modules
$M^{(g)}, W_1^{(g)}, W_2^{(g)}$  such that
\begin{equation}\label{iisso}(D_{g\m}\ot eBe)^g\!\cong M^{(g)}\ot  eBe,\,(D_{g\m}\ot A_1)^g\!\cong W_1^{(g)}
\ot  A_1,\,(D_{g\m}\ot A_2)^g\!\cong W_2^{(g)} \ot A_2, \end{equation}
for any $g\in G,$ as $R\ot_{R^\af}(eBe)^{\rm{op}}$-, $R\ot_{R^\af}(A_1)^{\rm{op}}$- and $R\ot_{R^\af}(A_2)^{\rm{op}}$-modules, respectively.

\vu

Let
$\tilde\af=(\tilde{D}_g, {\tilde\af}_g)_{g\in G}$ denote the partial action of $G$ on $R\ot_{R^\af}R,$ where
$$\tilde{D}_g=D_g\ot_{R^\af}D_g\, \,\text{and}\,\, {\tilde\af}_g\colon \tilde{D}_{g\m}\to \tilde{D}_g \,\,\text{is induced by}
\,\, x\ot_{R^\af}y\mt \afg(x)\ot_{R^\af}\afg(y). $$
Since $e_g=(1_g\ot_{R^\af}1_g)e$  satisfies
$$e_g(d\ot_{R^\af} 1_g- 1_g\ot_{R^\af}d)=e(d\ot_{R^\af} 1- 1\ot_{R^\af}d)(1_g\ot_{R^\af}1_g)=0,$$ for all $d\in D_g, g\in G,$
then $e_g$ is a  separability idempotent for the commutative $R^\af$-algebra $D_g.$ The  fact that $\tilde\af_g$ is a ring  isomorphism implies that
${\tilde\af}_g(e_{g\m})\in D_g\ot_{R^\af}D_g$ is another separability idempotent for $D_g$. Since separability idempotents
for commutative algebras are unique,
\begin{equation}\label{egg}{\tilde\af}_g(e_{g\m})=e_g,\,\,\, \text{for all}\,\,\,\, g\in G.\end{equation}

On the other hand,
$$(\tilde{D}_{g\m}\ot_{R^e}Be)^g=( ( D_{g\m}\ot_{R^\af}D_{g\m}) \ot_{R^e}Be)^g,$$
 is an  $R^e\ot_{R^\af}(eBe)^{op}$-module via the action  induced by
\begin{equation}\label{actiongg} [(r_1\ot_{R^\af}r_2)\ot_{R^\af}b]\bullet [(d_1\ot_{R^\af}d_2)\ot_{R^e}b')]=
\tilde\af_{g\m}[(r_1\ot_{R^\af}r_2)1_{\tilde g}](d_1\ot_{R^\af}d_2)\ot_{R^e}b'b,
\end{equation}
for all $r_1,r_2\in R, d_1,d_2\in D_{g\m}, \,b\in eBe,\,b'\in Be $ and $1_{\tilde g}=1_g\ot_{R^\af} 1_g, g\in G.$
\medskip

 Then, there is an $R^e\ot_{R^\af}(eBe)^{op}$-module isomorphism
\begin{equation*} (\tilde{D}_{g\m} \ot_{R^e}Be)^g\cong (\tilde{D}_{g\m} \ot_{R^e}B)^g e.
\end{equation*}

Notice that for any  $R^e\ot_{R^\af}(eBe)^{op}$-module $M,$ the  abelian group $(e\ot_{R^\af} 1_{eBe})M$ is an $ R\ot_{R^\af}(eBe)^{op}$-module via
$$(r\ot_{R^\af} ebe)\cdot((e\ot_{R^\af}1_{eBe})m)=((r\ot_{R^\af} 1_R)e\ot_{R^\af} ebe)m,$$
for all $m\in M,\, r\in R$ and $b\in B.$ In particular,  for $M=  (\tilde{D}_{g\m} \ot_{R^e}B)^g e$ we see that
\begin{align*}
(e\ot_{R^\af}1_{eBe})\bullet( \tilde{ D}_{g\m} \ot_{R^e}B)^g e
&=\tilde\af_{g\m}(e_g)( \tilde{D}_{g\m} \ot_{R^e}B)^g e\\
&=e_{g\m}(\tilde{ D}_{g\m} \ot_{R^e}B)^g e\\
&=e(\tilde{D}_{g\m} \ot_{R^e}B)^g e
\end{align*}
is an $R\ot_{R^\af}(eBe)^{\rm{op}}$-module via
\begin{align*}
(r\ot_{R^\af} ebe)\cdot(e(d_1\ot_{R^\af}d_2)\ot_{R^e} b'e) =
&((r\ot_{R^\af} 1_R)e\ot_{R^\af} ebe)\bullet (e(d_1\ot_{R^\af}d_2)\ot_{R^e} b'e)\\
\stackrel{\eqref{actiongg},\eqref{egg}} =&[\tilde\af_{g\m}((r\ot_{R^\af} 1)e_g)](d_1\ot_{R^\af}d_2)\ot_{R^e} b'ebe\\
=&(\af_{g\m}(r1_g)\ot_{R^\af} 1_R)e(d_1\ot_{R^\af}d_2)\ot_{R^e} b'ebe\\
=&\af_{g\m}(r1_g)e(d_1\ot_{R^\af}d_2)\ot_{R^e} b'ebe.
\end{align*}
Also $$[\tilde{D}_{g\m}\ot_{R^e}eBe]^g=[e\tilde{D}_{g\m}\ot_{R^e}Be]^g$$
 is an $ R\ot_{R^\af}(eBe)^{\rm{op}}$-module via
\begin{equation}\label{emod}(r\ot_{R^\af} ebe)  \blacktriangleright (e(d_1\ot_{R^\af}d_2)\ot_{R^e}b'e)=\af_{g\m}(r1_g)e(d_1\ot_{R^\af}d_2)\ot_{R^e}b'ebe.\end{equation}
We conclude that
\begin{equation}\label{isso} { e(\tilde{D}_{g\m} \ot_{R^e}B)^g e\cong [\tilde{D}_{g\m}\ot_{R^e}eBe]^g },\end{equation} as $ R\ot_{R^\af}(eBe)^{op}$-modules.
\smallskip

Moreover we have the next.

\begin{claim}\label{iso11} There is an $ R\ot_{R^\af}(eBe)^{\rm{op}}$-module isomorphism $$(D_{g\m}\ot eBe)^g\cong(\tilde{D}_{g\m}\ot_{R^e} eBe)^g,$$
where  $ R\ot_{R^\af}(eBe)^{\rm{op}}$ acts on $(\tilde{D}_{g\m}\ot_{R^e} eBe)^g$ as in {\rm(\ref{emod}).}
\end{claim}

Indeed, the map  $(D_{g\m}\ot eBe)^g \stackrel{\varsigma}{\to} (\tilde{D}_{g\m}\ot_{R^e} eBe)^g,$ determined by
$$d\ot ebe\mt (1_{g\m}\ot_{R^\af} d)\ot_{R^e} ebe$$
is a well defined $ { R  \ot }_{R^\af}(eBe)^{\rm{op}}$-module isomorphism whose inverse
 $$ (\tilde{D}_{g\m}\ot_{R^e} eBe)^g\stackrel{\varsigma^*}{\to} (D_{g\m}\ot eBe)^g,$$
is induced by
$$(d_1\ot_{R^\af} d_2)\ot_{R^e} ebe\mt d_1d_2\ot ebe.$$

In fact, since $e\in R^e$ is the separability idempotent for $R,$  we have $(d_1\ot_{R^\af} d_2)e=(1_{g\m}\ot_{R^\af}d_1 d_2)e.$ Hence,
\begin{align*}
\varsigma\varsigma^*((d_1\ot_{R^\af} d_2)\ot_{R^e} ebe)&=\varsigma(d_1d_2\ot ebe) =(1_{g\m}\ot_{R^\af} d_1d_2)\ot_{R^e} ebe \\
&=(1_{g\m}\ot_{R^\af} d_1d_2)e\ot_{R^e} ebe
=(d_1\ot_{R^\af} d_2)e\ot_{R^e} ebe\\
&=(d_1\ot_{R^\af} d_2)\ot_{R^e} ebe,
\end{align*}
and
$$\varsigma^*\varsigma(d_1\ot ebe)=\varsigma^*((1_{g\m}\ot_{R^\af } d_1)\ot_{R^e}  ebe)=d_1\ot ebe. $$

\vu

Now we prove that $\varsigma$ is ${R \ot}_{R^\af}(eBe)^{\rm{op}}$-linear. For
\begin{align*}
\varsigma((r\ot_{R^\af} ebe)\!\bullet\! (d\ot eb'e))&=\!\varsigma(\af_{g\m}(r1_g)d\ot eb'ebe)\!
=\!(1_{g\m}\!\ot_{R^\af}\af_{g\m}(r1_g)d)\ot_{R^e} eb'ebe,
\end{align*} and
\begin{align*}
(r\ot_{R^\af} ebe) \blacktriangleright \varsigma( d\ot eb'e)
&=(r\ot_{R^\af} ebe) \blacktriangleright ((1_{g\m}\ot_{R^\af } d)\ot_{R^e}  eb'e)\\
&=\af_{g\m}(r1_g)e(1_{g\m}\ot_{R^\af}d)\ot_{R^e}eb'ebe\\
&=(1_{g\m}\ot_{R^\af}\af_{g\m}(r1_g))e(1_{g\m}\ot_{R^\af}d)\ot_{R^e}eb'ebe\\
&=(1_{g\m}\ot_{R^\af}\af_{g\m}(r1_g))(1_{g\m}\ot_{R^\af}d)\ot_{R^e}eb'ebe\\
&=(1_{g\m}\ot_{R^\af}\af_{g\m}(r1_g)d) \ot_{R^e} eb'ebe,
\end{align*}
which ends the proof of the claim.

\vd

We still need the following.
\begin{claim}\label{1iso} There is an  $(R^e)\ot_{R^\af}(eBe)^{\rm{op}}$-module isomorphism
$$[\tilde{D}_{g\m}\ot_{R^e}(A_1\ot_{R^\af} A_2)]^g\cong (D_{g\m}\ot A_1)^g\ot_{R^\af} (D_{g\m}\ot A_2)^g,$$
where the action of $(R^e)\ot_{R^\af}(eBe)^{op}$ on $(D_{g\m}\ot A_1)^g\ot_{R^\af} (D_{g\m}\ot A_2)^g$ is induced by
\begin{align*}&[(r_1\ot_{R^\af}r_2)\ot_{R^\af} (x\ot_{R^\af}y)]\bullet [(d_1\ot a_1)\ot_{R^\af}(d_2\ot a_2)]\\
=&(\af_{g\m}(r_11_g)d_1\ot a_1x)\ot_{R^\af}(\af_{g\m}(r_21_g)d_2\ot a_2y),\end{align*}
for $r_1,r_2\in R, d_1,d_2\in D_{g\m}, x\ot_{R^\af}y\in eBe, a_1\in A_1$ and $a_2\in A_2.$
\end{claim}

\vu

Indeed, we have a well defined (additive) group homomorphism
$$\chi:(\tilde{D}_{g\m}\ot_{R^e} (A_1\ot_{R^\af} A_2))^g\to (D_{g\m}\ot A_1)^g\ot_{R^\af} (D_{g\m}\ot A_2)^g,$$ determined by
$$(d_i\ot_{R^\af} d'_i) \ot_{R^e} (a_j\ot_{R^\af}a'_j)\mt (d_i\ot a_j)\ot_{R^\af}(d'_i\ot a'_j), $$ which has
$$(d_i\ot a_i)\ot_{R^\af} (d_j\ot_{R^\af}a_j)\mt (d_i\ot_{R^\af} d_j)\ot_{R^e}(a_i\ot a_j), $$
 as an inverse. The fact that     $\chi$ is $(R^e)\ot_{R^\af}(eBe)^{op}$-linear is straightforward.

\vd

By (\ref{iisso}),  Claim \ref{iso11}, (\ref{isso})
and Claim \ref{1iso} we have the $R\ot_{R^\af} (eBe)^{\rm{op}}$-module isomorphisms
\begin{align*}
M^{(g)}\ot eBe&\cong (D_{g\m}\ot eBe)^g \cong[\tilde{D}_{g\m}\ot_{R^e} { eBe]^g}\\
&\cong e[\tilde{D}_{g\m}\ot_{R^e} { B]^g  }  e
\cong e[(D_{g\m}\ot { A_1)^g}  \ot_{R^\af} (D_{g\m}\ot {A_2)^g}  ]e\\
&\cong e[(W_1^{(g)}\ot A_1)\ot_{R^\af} (W_2^{(g)}\ot A_2)]e
\cong e\{W_1^{(g)}\ot [(A_1\ot_{R^\af} A_2)\ot W_2^{(g)}]\}e\\
&\cong e[(W_1^{(g)}\ot W_2^{(g)})\ot (A_1\ot_{R^\af} A_2)]e
\cong e[(W_1^{(g)}\ot W_2^{(g)})\ot (A_1\ot_{R^\af} A_2)e]\\
&\cong e[(A_1\ot_{R^\af} A_2)e \ot(W_1^{  (g) }\ot W_2^{(g)})]
\cong e(A_1\ot_{R^\af} A_2)e \ot(W_1^{ (g) } \ot W_2^{(g)})\\
&\cong (W_1^{(g)}\ot W_2^{(g)})\ot eBe.
\end{align*}

Finally Lemma \ref{cancel} implies $M^{(g)}\cong  W_1^{(g)} \ot W_2^{(g)},$ and
we conclude that $\varphi_5$ is a group homomorphism.\cua

\vd

\section{Two partial representations   $G\to{\bf PicS}_{R^\af}(R)$ and the homomorphism
$H^1(G,\af^*,{\bf PicS}(R)) \stackrel{\f_6}\to  H^3(G,\af, R)$}\label{phi6}

\vu

 For reader's convenience we recall from \cite{DEP} the concept of a partial representation.

\begin{defi}\label{defn-par-repr} A (unital) partial representation of $G$ into an algebra (or, more ge\-ne\-ral\-ly, a monoid) $S$ is a map $\Phi: G \to S$ which satisfies the following properties, for all $g,h\in G,$

\vspace*{2mm}
\noindent (i) $\Phi (g\m) \Phi (g) \Phi (h) = \Phi (g\m) \Phi (g h),$

\vspace*{2mm}
\noindent (ii) $\Phi (g ) \Phi (h) \Phi (h\m) = \Phi (g h) \Phi ( h \m ),$

\vspace*{2mm}
\noindent (iii) $\Phi (1_G )  = 1_S.$
\end{defi} 

For any $g\in G$ denote by $_{g}{(D_{g\m})}_I$ the $R$-$R$-bimodule $D_{g\m}$ regarded as an  $R$-$R$-bimodule with new action $*$ given by
$$r*d=\af_{g\m}(r1_{g})d,\,\,\,\,\text{and} \,\,\, d*r=dr, \,\,\, \text{for any}\,\,\, r\in R,\,d\in D_{g\m},$$ (analogously
we define ${}_{I}{(D_{g})}_{g\m};$ notice that $D_g={}_{I}{(D_{g})}_I$ as $R$-$R$-bimodules). Then using (iii) of Lemma \ref{gaction},
we get that $ [_{g}{(D_{g\m})}_I]\in {\bf PicS}_{R^\af}(R).$ We set
\begin{equation}\label{fi0}\Phi_0\colon G\ni g\mt [_{g}{(D_{g\m})}_I]\in {\bf PicS}_{R^\af}(R).\end{equation}
Some useful properties of $\Phi_0$ are given in the next.

\begin{prop}\label{pr1} Let $\Phi_0$ be as in {\rm(\ref{fi0})}. Then,
\begin{itemize}
\item $\Phi_0$ is a partial representation of $G$ in ${\bf PicS}_{R^\af}(R)$ with
\begin{equation}\label{dg}\Phi_0(g)\Phi_0(g\m)=[D_g],\end{equation} for all $g\in G,$

\vu

\item $\Phi_0(g)[D_{h}]=[D_{gh}]\Phi_0(g),$ for any $g,h\in G.$ In particular, \begin{equation}\label{fiab}\Phi_0(g)[D_{g\m}]=\Phi_0(g)=[D_g]\Phi_0(g),\end{equation}

\vu

\item for any $g\in G$ and  $[P]\in X_{g\m}$ there is an $R-R$-bimodule isomorphism \begin{equation}\label{isoa}P_g\cong \Phi_0(g)\ot P\ot \Phi_0(g\m).\end{equation}
\end{itemize}

\vu

\end{prop}
\p It is clear that $\Phi_0(1)=[R].$ Now let $g,h\in G$ and consider the map
 $\phi_{g,h}\colon _{g}{(D_{g\m})}_I \ot\, _{h}{(D_{h\m})}_I \ot\,  _{h\m}{(D_{h})}_I \to   _{gh}\!\!(D_{(gh)\m})_I
 \ot\,  _{h\m}{(D_{h})}_I,  $ defined by $$\phi_{g,h}(a_{g\m}\ot b_{h\m} \ot c_h)= \af_{h\m}(a_{g\m}1_h)b_{h\m}\ot c_h.$$
 
 Then,  $\phi_{g,h}$   is $R$-$R$-linear because for $r_1,r_2\in R$ we have
\begin{align*}\phi_{g,h}(r_1 *a_{g\m}\ot b_{h\m} \ot c_h \ast r_2)&=\phi_{g,h}(\af_{g\m}(r_11_{g})a_{g\m}\ot b_{h\m} \ot c_hr_2)\\
&=\af_{h\m}(\af_{g\m}(r_11_{g})a_{g\m}1_h)b_{h\m}\ot c_hr_2\\
&=\af_{(gh)\m}(r_11_{gh})\af_{h\m}(a_{g\m}1_h)b_{h\m}\ot c_hr_2\\
&=r_1*\phi_{g,h}(a_{g\m}\ot b_{h\m} \ot c_h) \ast r_2.\end{align*}
Moreover,  the map   $ _{gh}{(D_{(gh)\m})}_I \ot\,  _{h\m}{(D_{h})}_I\to{}_{g}{(D_{g\m})}_I \ot{}_{h}{(D_{h\m})}_I \ot{}_{h\m}{(D_{h})}_I$
induced by $$x_{(gh)\m}\ot x_h\to \af_h(x_{(gh)\m}1_{h\m})\ot 1_{h\m}\ot x_h,$$ for all $x_{(gh)\m}\in\, _{gh}{(D_{(gh)\m})}_I,\, x_h\in _{h\m}(D_{h})_I $
 is the inverse of $\phi_{g,h}.$ This  yields that $\Phi_0(g)\Phi_0(h)\Phi_0(h\m)=\Phi_0(gh)\Phi_0(h\m).$

\vd

In a similar way,  the map
$ _{g\m}{(D_{g})}_I \ot\, _{gh}\!(D_{(gh)\m})_I\to  _{g\m}\!\!(D_{g})_I \ot\, _{g}{(D_{g\m})}_I \ot _{h}{(D_{h\m})}_I,  $ such that
 $$a_{g}\ot b_{(gh)\m} \mt a_{g} \ot\af_h(b_{(gh)\m}1_{h\m}) \ot 1_{h\m},$$ is an $R$-$R$-bimodule isomorphism with inverse
$$ x_g \otimes y_{g\m} \ot z_{h\m} \mapsto x_g \ot \alpha _{h\m} (y_{g\m } 1_h) z_{h\m},$$  and we obtain
 $\Phi_0(g\m)\Phi_0(gh)=\Phi_0(g\m)\Phi_0(g)\Phi_0(h).$

\vd

To prove (\ref{dg}) one can check that  the map determined by  $${}_{g}{(D_{g\m})}_I\, \ot\, _{g\m}{(D_{g})}_I \ni a_{g\m}\ot b_{g} \mt \af_g(a_{g\m}) b_{g}\in D_{g}$$
is a well defined $R$-$R$-bimodule isomorphism whose inverse is
$D_g \ni d\to 1_{g\m}\ot d\in\, _{g}{(D_{g\m})}_I\, \ot {}_{g\m}{(D_{g})}_I.$

\vd

The second item follows from the first and (2), (3) of \cite{DEP}. To check the last item consider the map
$$P_g\ni p\stackrel{\nu}{\to} 1_{g\m}\ot p\ot 1_g\in{}_{g}{(D_{g\m})}_I \ot \, P \ot{}
_{g\m}{(D_{g})}_I.$$ Then, for $r_1,r_2\in R$ we have
\begin{align*}\nu(r_1\bullet p\bullet r_2) &=1_{g\m}\ot r_1\bullet p\bullet r_2\ot 1_g=1_{g\m} \ot \af_{g\m}(r_11_g)p\af_{g\m}(r_21_g)\ot 1_g\\
&=\af_{g\m}(r_11_g)\ot p\ot r_21_g= r_1\ast 1_{g\m}\ot p\ot  1_g \ast r_2=r_1 \ast \nu(p) \ast r_2.\end{align*}
Hence, $\nu$  is an $R$-$R$-bimodule isomorphism with inverse  induced by $$a_{g\m}\ot p\ot b_g =
a_{g\m}\ot p\ot \af_{g\m}( b_g) \ast 1_g = 1_{g\m} \ot a_{g\m}p\af_{g\m}(b_g) \ot 1_g \mt a_{g\m}p\af_{g\m}(b_g),$$
for all $a_{g\m}\ot p\ot b_g\in\!  _{g}{(D_{g\m})}_I \ot\, P \ot \,  _{g\m}{(D_{g})}_I.$\cua

\vt

Now we construct another partial representation of  $G$ in  ${\bf PicS}_{R^\af}(R).$
\begin{lema} \label{pr2}For any $f\in Z^1(G,\af^*, {\bf PicS}(R))$ set $\Phi_f=f\Phi_0\colon G\to  {\bf PicS}_{R^\af}(R),$
that is $\Phi_f(g)=f(g)\Phi_0(g),$ for any $g\in G.$ Then,
\begin{itemize}
\item $\Phi_f$ is a partial representation,

\vu

\item $\Phi_f(g)\Phi_f(g\m)=[D_g],$ for all $g\in G.$
\end{itemize}

Moreover, 
 writing 
$\Phi_f(g)=[J_g],$  we have that  $D_g\cong {\rm End}_{D_g}(J_g),$ as $R$-
and   $D_g$-algebras,  for any $g\in G.$

\end{lema}

\p First of all we have $\Phi_f(1)=[R].$ Now, let $g,h\in G.$ Then,
\begin{align*}\Phi_f(g\m)\Phi_f(gh)&=f(g\m)\Phi_0(g\m)f(gh)\Phi_0(gh)\\
&\stackrel
{(\ref{fiab})}=f(g\m)\Phi_0(g\m)[D_g]f(gh)\Phi_0(gh)\\
&=f(g\m)\Phi_0(g\m)f(g)\af^*_g(f(h)[D_{g\m}])\Phi_0(gh)\\
&\stackrel{(\ref{isoa})}{=}f(g\m)\Phi_0(g\m)f(g)\Phi_0(g)f(h)[D_{g\m}]\Phi_0(g\m)\Phi_0(gh) \\
&\stackrel{(\ref{fiab})}=f(g\m)\Phi_0(g\m)f(g)\Phi_0(g)f(h)\Phi_0(g\m)\Phi_0(gh)\\
&=f(g\m)\Phi_0(g\m)f(g)\Phi_0(g)f(h)\Phi_0(g\m)\Phi_0(g)\Phi_0(h)\\
&\stackrel{(\ref{dg})}=f(g\m)\Phi_0(g\m)f(g)\Phi_0(g)f(h)[D_{g\m}]\Phi_0(h)\\
&\stackrel{(\ref{fiab})}=f(g\m)\Phi_0(g\m)f(g)\Phi_0(g)f(h)\Phi_0(h)\\
&=\Phi_f(g\m)\Phi_f(g)\Phi_f(h).\end{align*}

Analogously, it can be shown that $\Phi_f(gh)\Phi_f(h\m)=\Phi_f(g)\Phi_f(h)\Phi_f(h\m).$ Indeed, since $f\in  C^1 (G,\af^*,{\bf PicS}(R))$
we have $[D_{h\m}]f(h\m)=f(h\m),$ and by the second item of Proposition \ref{pr1} we obtain $\Phi_0(gh)[D_{h\m}]=[D_g]\Phi_0(gh).$ Thus,
\begin{align*}\Phi_f(gh)\Phi_f(h\m)&=f(gh)\Phi_0(gh)f(h\m)\Phi_0(h\m)\\
&=f(gh)[D_g]\Phi_0(gh)f(h\m)\Phi_0(h\m)\\
&\stackrel{(\ref{isoa})}{=}f(g)\Phi_0(g)f(h)[D_{g\m}]\Phi_0(g\m)\Phi_0(gh)f(h\m)\Phi_0(h\m)\\
&\stackrel{(\ref{fiab})}=f(g)\Phi_0(g)f(h)\Phi_0(g\m)\Phi_0(gh)f(h\m)\Phi_0(h\m)\\
&\stackrel{(\ref{dg})}=f(g)\Phi_0(g)f(h)[D_{g\m}]\Phi_0(h)f(h\m)\Phi_0(h\m)\\
&=f(g)\Phi_0(g)f(h)\Phi_0(h)f(h\m)\Phi_0(h\m)\\
&=\Phi_f(g)\Phi_f(h)\Phi_f(h\m).
\end{align*}

With respect to the second item we have $\Phi_f(g)\Phi_f(g\m)=f(g)\Phi_0(g)f(g\m)\Phi_0(g\m)\stackrel{(\ref{isoa})}{=}f(g)\af^*_g(f(g\m))=f(gg\m)[D_g]=[R][D_g]=[D_g].$

\vd

Finally, 
by definition
$[J_g]= f(g) \Phi_0(g),$ and  $J_g\cong D_g\ot J_g$  as left $R$- and $D_g$-modules,  for all $g\in G$.
Therefore, after identifying $D_g={\rm End}_{D_g}(D_g)$,  the $R$-algebra epimorphism $ R\to {\rm End}_R(J_g)$
induces an $R$- and $D_g$-algebra epimorphism $\xi\colon D_g\to {\rm End}_{D_g}(J_g),$ thanks to Proposition \ref{ht}.
Via localization we will check that $\xi$ is injective.

\vu

 If for any $g\in G$ the ring $D_g$ is semi-local,  { then  ${\bf Pic}(D_g)$ is trivial (see, for example,  \cite[Ex. 2.22 (D)]{L}), and 
by Remark \ref{unxg}   we have that} 
 $f(g)\cong D_g,$ as $ D_g$-modules  as well as $ R $-modules. Then
\begin{equation}\label{anj}J_g\cong \, _g (D_{g\m})_I\,\,\,\text{as  $R$-$R$-bimodules.}\end{equation}
  Moreover, after localizing by a prime ideal of $R^\af , $ we obtain $R^\af$-module isomorphisms $D_{g\m}\cong (R^\af)^m\cong D_g$ for some $m\in\N,$  thanks to the facts that  the maps $\alpha_g$ are isomorphisms of $R^\af$-modules and localization is an exact functor. Since in this case any $D_g$ is semi-local, \eqref{anj} implies that the map 
 $\xi \colon (R^\af)^m\to
{\rm End}_{(R^\af)^m}(_g(R^\af)^m)$ is an epimorphism of $R^\af$-algebras. On the other hand, the left $R^\af$-modules $(R^\af)^m$ and $_g(R^\af)^m$ are isomorphic
via $r \mapsto {\alpha}_{g\m} (r),$ and we have ${\rm End}_{(R^\af)^m}(_g(R^\af)^m) \cong (R^\af)^m,$ and $\xi $ must be an isomorphism of $R^\af$-modules. Finally, since $\xi $ is $R$-linear, then it is an isomorphism of $R$-
and   $D_g$-algebras,  for any $g\in G.$
\cua

\vd

In what follows we shall write  $\Phi_f (g) =[J_g].$

\vu

\begin{remark}\label{jgfiel} Localizing by ideals in ${\rm Spec}(R^\af)$ and using   \eqref{anj} we see that    $J_g$ is a faithful $D_g$-module.
Then if we ignore the right $R$-module structure of  $J_g$ and use  the last item of Lemma~\ref{pr2} we obtain
that $[J_g]\in {\bf Pic}(D_g),$ for any $f\in Z^1(G,\af^*, {\bf PicS}(R))$ and $g\in G.$ In particular,
the  map $m_{D_g}\colon D_g\to  {\rm End}_{D_g}(J_g)$ given by left multiplication  is a $D_g$-algebra isomorphism.
\end{remark}
\vu

\begin{remark}\label{jgdg} Let $f$ be an element of $ Z^1(G,\af^*, {\bf PicS}(R))$ and  write  $f(g)=[ M_g ].$
Then $J_g= M_g \ot \, _{g}{(D_{g\m})}_I$ and
\begin{equation}\label{jgtwits}x_gr=\af_g(r1_{g\m})x_g,\,\,\text{for any}\,\, x_g\in J_g,\,\, r\in R,\end{equation} and the map
\begin{equation}\label{jgdgi} J_g\ni x\to x\ot 1_{g\m}\in J_g\ot D_{g\m}
\end{equation}
is an $R$-$R$-bimodule isomorphism. Furthermore, if $D_g$ is a semi-local for any $g\in G,$ then by \eqref{anj} 
 there is an $R$-$R$-bimodule isomorphism $\gamma_g\colon _{g}{(D_{g\m})}_I\ni d\to \gamma_g(d)\in J_g .$ Therefore, setting
$u_g=\gamma_g(1_{g\m}),$ we have, for any $x\in J_g,$ that  $$x=\gamma_g(d)=\gamma_g(\af_g(d)* 1_{g\m})=\af_g(d)u_g,$$ for some $d\in D_{g\m}.$
 We conclude that $J_g=D_gu_g,$ $u_g$ is a free  generator of $J_g$ over $D_g$ and
 $u_gr=\af_g(r1_{g\m})u_g,$ for all $ r\in R,\,g\in G,$ in view of (\ref{jgtwits}).
\end{remark}

\vu

We know  from Proposition~\ref{pr1} and Lemma~\ref{pr2} that there are   $R$-$R$-bimodule isomorphisms
${}_{h\m}(D_{h})_{I}\ot {}_{g\m}(D_{g})_I\cong D_{h\m}\ot \,{}_{h\m g\m}(D_{g h})_{I}\,\, \,\,\text{and}\,\,\,\, J_g\ot \, D_h \cong D_{gh}\ot J_g,$
 for all $g,h\in G.$ For further reference, we shall construct these isomorphisms explicitly.
\begin{lema}\label{iisoo} The map  $\varrho\colon {}_{h\m}(D_{h})_{I}\ot {}_{g\m}(D_{g})_I\to D_{h\m}\ot \,{}_{h\m g\m}(D_{g h})_{I},$ induced by
 \begin{equation*}\label{isovr}x_h\ot y_g\mt 1_{h\m}\ot \af_g(x_h1_{g\m})y_g,
\end{equation*}
is an $R$-$R$-bimodule isomorphism, and $\varrho\m\colon  D_{h\m}\ot\,{}_{h\m g\m}(D_{g h})_{I}\to {}_{h\m}(D_{h})_{I}\ot {}_{g\m}(D_{g})_I$
is determined by $$\label{isovri}x_{h\m}\ot  y_{gh}\mt\af_{g\m}(y_{gh}1_g)\af_h(x_{h\m})\ot 1_{g},$$
for $g,h\in G.$
\end{lema}

\p Indeed,  $\varrho$ is well defined, and using \eqref{prodp}  one can show that $\varrho$ is an $R$-$R$-bimodule homomorphism. Moreover,
\begin{align*} x_h\ot y_g&\stackrel{\varrho}\mt 1_{h\m}\ot \af_g(x_h1_{g\m})y_g \stackrel{\varrho\m}{\mt}\af_{g\m}(\af_g(x_h1_{g\m})y_g)1_h\ot 1_g
=x_h\af_{g\m}(y_g)\ot 1_g
=x_h\ot y_g.\end{align*}
On the other hand,
\begin{align*}x_{h\m}\ot y_{gh}&\stackrel{\varrho\m}\mt \af_{g\m}(y_{gh}1_g)\af_h(x_{h\m})\ot 1_{g}
\stackrel{\varrho}\mt1_{h\m}\ot \af_g(\af_{g\m}(y_{gh}1_g)\af_h(x_{h\m}))\\
&=1_{h\m}\ot  y_{gh}\af_g(\af_h(x_{h\m})1_{g\m})
\stackrel{\eqref{prodp}}=1_{h\m}\ot  y_{gh}\af_{gh}(x_{h\m}1_{({gh})\m})1_{g}\\
&=1_{h\m}\ot y_{gh}\af_{gh}(x_{h\m}1_{(gh)\m})\af_{gh}(1_{h\m}1_{(gh)\m})
=1_{h\m}\ot y_{gh}\af_{gh}(x_{h\m}1_{(gh)\m})\\
&=1_{h\m}\ot  (x_{h\m}*y_{gh})
=x_{h\m}\ot y_{gh},\end{align*}
as desired.\cua
\begin{lema}\label{kapa} The map $ J_g\ot  D_h \stackrel{\kappa_{g,h}}\to D_{gh}\ot  J_g$ induced by
\begin{equation}\label{kappa}a_{g}\ot \, b_h\mt\af_g(b_h1_{g\m})\ot a_{g},\end{equation}
for any $g,h\in G,$ is an $R$-$R$-bimodule isomorphism.
\end{lema}
\p First,  $\kappa_{g,h}$ is well defined by (\ref{jgtwits}). Notice that $\kappa_{g,h}$ is bijective with inverse
$\iota_{g,h}\colon D_{gh}\ot  J_g\to J_g\ot \, D_h, $  $a_{gh}\ot b_g\mt b_g\ot \af_{g\m}(a_{gh}1_g),$ for all $a_{gh}\in D_{gh}$ and $b_g \in J_g.$  Indeed,
$$\iota_{g,h}\circ \kappa_{g,h}(a_{g}\ot \, b_h)=\iota_{g,h}(\af_g(b_h1_{g\m})\ot a_{g})=a_{g}\ot b_h1_{g\m}\stackrel{\eqref{jgtwits}}=1_ga_{g}\ot b_h =a_{g}\ot b_h.$$
In addition,
$$ \kappa_{g,h}\circ \iota_{g,h}(a_{gh}\ot  b_g)= \kappa_{g,h}( b_g\ot \af_{g\m}(a_{gh}1_g))=a_{gh}1_g\ot  b_g= a_{gh} \ot  1_g b_g= a_{gh}\ot b_g.$$

Finally,  to prove that $\kappa_{g,h}$ is $R$-$R$-linear take $r_1, r_2\in R.$ Then,
\begin{align*}
\kappa_{g,h}(r_1\cdot a_{g}\ot \, b_hr_2)&=\af_g(b_hr_21_{g\m})\ot r_1\cdot a_{g}=\af_g(b_h1_{g\m})\af_g(r_21_{g\m})\ot r_1\cdot a_{g}\\
&=\af_g(b_h1_{g\m})r_1\ot \af_g(r_21_{g\m})\cdot a_{g}=r_1\af_g(b_h1_{g\m})\ot  a_{g}r_2.\end{align*}
This completes the proof.\cua

\subsection{The map $H^1(G,\af^*,{\bf PicS}(R)) \stackrel{\f_6}\to  H^3(G,\af, R)$} Let $f\in Z^1(G,\af^*, {\bf PicS}(R))$
and  $\Phi_f=f\Phi_0.$ Write $ \Phi _f (g)=[J_g].$  By Lemma \ref{pr2} the map  $\Phi_f$ is a partial homomorphism such
that $\Phi_f(g)\Phi_f(g\m)=[D_g].$ Then, there is a family of $R$-$R$-bimodule isomorphisms $$\{\chi  _{g,h}\colon J_g\ot J_h\to D_g\ot J_{gh}\}_{g,h\in G}.$$

Consider the following diagram
\begin{equation} \label{diag1}
\xymatrix{J_g\ot J_h\ot  J_l\ar[d]^{ \chi_{g,h}\ot{\rm id}_l}\ar[r]^{{\rm id}_g\ot \chi_{h,l}}& J_g\ot  D_{h}\ot J_{hl}
\ar[r]^{\kappa_{g,h}\ot {\rm id }_{hl}}&\,\,\,D_{gh}\ot J_g \ot J_{hl}\,\,\ar[r]^{{\rm id}_{D_{gh}}\ot \chi_{g,hl}} &\,\,\,\,D_{gh}
\ot  D_g\ot  J_{ghl}\ar[d]^{{\tau}_{gh,g}\ot {\rm id}_{ghl}} \\
D_g\ot J_{gh}\ot J_l\ar[rrr]^{{\rm id}_{D_g}\ot \chi_{gh,l}}&&&D_{g}\ot D_{gh}\ot J_{ghl}},
\end{equation}
for any $g,h,l\in G,$ where $\kappa_{g,h}$ is from Lemma \ref{kapa} and ${\tau}_{gh,g}$ is the twisting  $u\ot v\mt v\ot u.$

\vd

We use diagram \eqref{diag1} to  construct  a cocycle in $Z^3(G,\af,R).$  Let $\tilde\om(g,h,l)$ be the map obtained making
a counterclockwise loop in \eqref{diag1}:
$$({\rm id}_{D_g}\ot\chi_{gh,l})\circ(\chi_{g,h}\ot{\rm id}_l)\circ({\rm id}_g\ot \chi_{h,l})\m\circ(\kappa_{g,h}\ot
{\rm id }_{hl})\m\circ({\rm id}_{D_{gh}}\ot \chi_{g,hl})\m\circ({\tau}_{gh,g}\ot {\rm id}_{ghl})\m,$$
for all $g,h,l\in G.$ Evidently,  $\tilde\om(g,h,l)$ is a left $R$-linear  automorphism of  $D_g\ot D_{gh}\ot  J_{ghl}.$

\vu

 Moreover, from the fact that
\begin{equation}\label{acom}(a_g\ot b_{gh}\ot  c_{ghl})\cdot(t_g\ot u_{gh}\ot v_{ghl})=a_gt_g\ot b_{gh}u_{gh}\ot
c_{ghl}v_{ghl}=a_gb_{gh}c_{ghl}t_g\ot u_{gh}\ot v_{ghl},\end{equation}
for all $a_g,t_g\in D_g, b_{gh},u_{gh}\in D_{gh},c_{ghl}\in D_{ghl},$ $v_{ghl}\in J_{ghl},$ and $g,h,l\in G,$ we conclude that
$\tilde\om(g,h,l)$ is an invertible element of $${\rm End}_{D_g\ot D_{gh}\ot D_{ghl}}(D_g\ot D_{gh}\ot
J_{ghl})\cong D_g\ot D_{gh}\ot {\rm End}_{D_{ghl}}( J_{ghl})\cong D_g\ot D_{gh}\ot  D_{ghl},$$ where the  last ring isomorphism
follows from  Remark \ref{jgfiel}.  Thus there is a unique in\-ver\-ti\-ble element $\om_1(g,h,l)\in \U(D_g\ot D_{gh}\ot D_{ghl})$ such that
 $\tilde\om(g,h,l)(z)=\om_1(g,h,l)z,$ for all $z\in D_g\ot D_{gh}\ot J_{ghl},$  and it follows from (\ref{acom})
 that there is a unique  $\om(g,h,l)\in \U(D_gD_{gh}D_{ghl})$ satisfying $$\tilde\om(g,h,l)z=\om(g,h,l)z,\,\,\, g,h,l\in G,\,\,\, z\in D_g\ot  D_{gh}\ot J_{ghl}. $$

We shall check that $\om\in Z^3(G,\af, R),$ or equivalently
\begin{equation}\label{cobordo}(\de^3\om)(g,h,l,t)=1_g1_{gh}1_{ghl}1_{ghlt},\,\,\text{ for all}\,\, g,h,l,t\in G,\end{equation}
where $\de^3$ is the coboundary operator given by (\ref{pcob}).

\vu

 Since $(\de^3\om)(g,h,l,t)$ and $1_g1_{gh}1_{ghl}1_{ghlt}$ belong to  the  $R^\af$-module $R,$  equality  \eqref{cobordo}  holds if and only if  for every $\mfp\in {\rm Spec}(R^\af)$
 the image of $( \de^3 \om )(g,h,l,t) $ in $(D_gD_{gh}D_{ghl}D_{ghlt})_\mfp$ is $(1_g1_{gh}1_{ghl}1_{ghlt})_\mfp.$
 But if $D_g$ is semi-local for any $g\in G$,  Remark \ref{jgdg} implies  $J_g=D_gu_g,$ and it follows that
 $$\chi_{g,h}(a_gu_g\ot b_h u_h)=a_g\af_g(b_h1_{g\m})\chi_{g,h}(u_g\ot u_h)=a_g\af_g(b_h1_{g\m})\tilde\rho(g,h) (1_g\ot  u_{gh} ),$$
 with $\tilde\rho(g,h)\in \U(D_g\ot D_{gh}).$ One can write $\tilde\rho(g,h)=1_g\ot \rho(g,h),$ where $\rho(g,h)$ belongs to $\U(D_gD_{gh}).$
 In particular, we have a map $\rho\in C^2(G,\af, R)$ and
\begin{equation*}\label{coci}\chi_{g,h}(a_gu_g\ot b_h { u_h})=a_g\af_g(b_h1_{g\m})\rho(g,h)1_g\ot u_{gh}.\end{equation*}
 We conclude that $\chi\m_{g,h}(1_g\ot u_{gh})=\rho(g,h)\m u_g\ot u_h,$ where  $\rho(g,h)\m$ is the inverse of $\rho(g,h)$ in $D_gD_{gh}.$
Now we apply  $\om(g,h,l)$ to $1_{g}\ot 1_{gh}\ot u_{ghl}.$
\begin{align*} 1_{g}\ot 1_{gh}\ot  u_{ghl}&\mt 1_{gh}\ot 1_{g}\ot u_{ghl}
\mt 1_{gh}\ot \rho\m(g,hl)u_g\ot u_{hl}
\mt \rho(g,hl)\m u_g\ot 1_h1_{g\m}\ot u_{hl}\\
&\stackrel{\eqref{jgtwits}}=\rho(g,hl)\m u_g\ot 1_h\ot u_{hl}
\mt\rho(g,hl)\m \af_g(\rho(h,l)\m1_{g\m})u_g\ot u_h\ot u_{l}\\
& \mt\rho(g,hl)\m \af_g(\rho(h,l)\m1_{g\m})\rho(g,h)1_g\ot  u_{gh}\ot u_{l}\\
& \mt\rho(g,hl)\m \af_g(\rho(h,l)\m1_{g\m})\rho(g,h)\rho(gh,l)1_g\ot 1_{gh}\ot u_{ghl}.
\end{align*}
 Thus,  \begin{equation*}\label{3cocl}\om(g,h,l)(1_g\ot 1_{gh}\ot u_{ghl})=(\de^2 {  \rho \m } )(g,h,l)(1_g\ot 1_{gh}\ot u_{ghl}), \, g,h,l\in G.\end{equation*}
 Hence, $\om(g,h,l)=(\de^2   {  \rho \m }  )(g,h,l)$ and it follows from Proposition \ref{pcobh} that $$(\de^3\om)(g,h,l,t)=1_g1_{gh}1_{ghl}1_{ghlt}.$$
 This yields  $\om\in Z^3(G,\af,R).$
\begin{claim} The map $\varphi_6 \colon H^1(G,\af^*,{\bf PicS}(R))\ni {\rm cls}(f)\to {\rm cls}(\om) \in H^3(G,\af,R)$ is well defined.
\end{claim}
\proof If one takes another family  $\{\chi'_{g,h}\colon J_g\ot  J_h\to D_g\ot J_{gh}\}_{g,h\in G}$ of $R$-$R$-bimodule isomorphisms,
the map $\chi'_{g,h}\circ \chi\m_{g,h}$ is an invertible element of ${\rm End}_{D_g\ot D_{gh}}(D_g\ot J_{gh}),$ and there exists
$\sigma(g,h)\in \U(D_gD_{gh})$ such that  $\chi'_{g,h}\circ \chi\m_{g,h}(z)=\sigma(g,h)z,\,\, \text{
 for all}\,\, z\in D_g\ot  J_{gh}.$
Thus, $\sigma\in C^2(G,\af, R),\, \chi_{g,h}=\sigma(g,h)\chi'_{g,h},$ and
setting $\tilde\om'(g,h,l)=({\rm id}_{D_g}\ot \chi'_{gh,l})\circ(\chi'_{g,h}\ot{\rm id}_l)\circ({\rm id}_g\ot \chi'_{h,l})\m\circ(\kappa_{g,h}\ot
{\rm id }_{hl})\m\circ({\rm id}_{D_{gh}}\ot \chi'_{g,hl})\m\circ({\tau}_{gh,g}\ot {\rm id}_{ghl})\m,$   we see that after localizaing by  ideals in ${\rm Spec}(R^\af)$
that \begin{equation*}\label{prodfam}\om'(g,h,l)=(\de^2\sigma\m)(g,h,l)\om(g,h,l).\end{equation*} This implies that ${\rm cls}(\om)={\rm cls}(\om')$ in $H^3(G,\af,R).$

\vu

On the other hand taking another representative $J'_g\in [J_g],$ for any  $g\in G,$ we have  families of $R$-$R$-bimodule isomorphisms
$\{\chi'_{g,h}\colon J'_g\ot  J'_h\to D_g\ot J'_{gh}\}_{g,h\in G}$ and $\{\zeta_g\colon J_g\to J'_g\}_{g\in G}.$

\vu

 Let
$\chi''_{g,h}= ({\rm id}_{D_g}\ot \zeta_{gh})\circ\chi_{g,h}\circ (\zeta\m_g\ot \zeta\m_h),\,\, g,h\in G.$ Thus,
if $\om'$ and $\om''$ are the corresponding cocycles in $Z^3(G,\af,R)$ induced by the families
$\{\chi'_{g,h}\colon J'_g\ot  J'_h\to D_g\ot J'_{gh}\}_{g,h\in G}$ and $\{\chi''_{g,h}\colon J'_g\ot J'_h\to D_g\ot J'_{gh}\}_{g,h\in G}$
respectively,  by the above we have ${\rm cls}(\om')={\rm cls}(\om'')$ in $H^3(G,\af,R).$ We shall prove that  $\om=\om''.$

\vu

By localization we may assume that  each $D_g,$ $g\in G,$ is a semi-local ring. Then, by Remark \ref{jgdg} there is $u_g\in J_g$  such that $J_g=D_gu_g$ and $J'_g=D_gu'_g,$
where $\zeta_g(u_g)=u'_g.$ Hence, the equality  ${ \chi _{g,h}  } (au_g\ot bu_h)=a\af_g(b1_{g\m})\rho(g,h)1_g\ot u_{gh},$ for all $g,h\in G,$ and the
definition of $\chi''_{g,h}$ imply   $$\chi''_{g,h}(au'_g\ot  bu'_h)=a\af_g(b1_{g\m})\rho(g,h)1_g\ot u'_{gh},$$ and using the construction
of $\om$ and $\om''$ we get $\om=\om''.$

\vu

Finally, if ${\rm cls}(f)={\rm cls}(f')\in H^1(G,\af^*,{\bf PicS}(R)),$  there exists $f_0\in B^1(G,\af^*,{\bf PicS}(R)) $
such that $f'=f_0f,$  and $[P]\in{\bf Pic}(R)$ with $f_0(g)=[P]\af^*_g([P^*][D_{g\m}]),$ for all $g\in G.$ Hence,
\begin{align*}\Phi_{f'}(g)= f'(g)\Phi_0(g)&=f_0(g)f(g)\Phi_0(g)\\
&=[P]\af^*_g([P^*][D_{g\m}])f(g)\Phi_0(g)\\
&\stackrel{(\ref{fiab}),(\ref{isoa})}=[P]\Phi_0(g)[P^*]\underbrace{\Phi_0(g\m)f(g)\Phi_0(g)}_{\in\, {\rm PicS}(R)}\\
&=[P]\Phi_0(g)\Phi_0(g\m)f(g)\Phi_0(g)[P^*]\\
&\stackrel{(\ref{dg})}=[P]f(g)\Phi_0(g)[P^*]\\
&=[P]\Phi_f(g)[P^*].
\end{align*}

Set $\Phi_f(g)=[J_g].$ Let $\{\chi _{g,h}\colon J_g\ot  J_h\to D_g\ot J_{gh}\}_{g,h\in G}$ be a family of $R$-$R$-bimodule
isomorphisms which come from $\Phi_f $, and $\om\in Z^3(G,\af,R)$ determined by  the $ \chi _{g,h}.$  Identifying
$(P\ot  J_g\ot  P^*)\ot (P\ot  J_h\ot  P^*)\cong P\ot  J_g\ot  J_h\ot  P^*$, we choose the family of $R$-$R$-bimodule isomorphisms
$$\{{\chi}'_{g,h}\colon P\ot J_g\ot  J_h\ot P^*\to P\ot  D_g\ot  J_{gh}\ot  P^*\}_{g,h\in G},$$ where
 ${\chi}'_{g,h}={\rm id}_P\ot \chi _{g,h}\ot{\rm id}_{P^*}\colon P\ot  J_g\ot  J_h\ot P^*\to P\ot  D_g\ot  J_{gh}\ot  P^*,$ for all $g,h\in G.$
 The isomorphisms  $\{{\chi}'_{g,h}\}$ correspond to  $\Phi_{f'}(g).$  Therefore, if $\om'\in Z^3(G,\af,R)$ is induced by the family $\{{\chi}'_{g,h}\}_{g,h\in G},$ we get
$$\om'(g,h,l)\!=\!{\rm id}_P\ot \om(g,h,l)\ot {\rm id}_{P^*}\!\!\in\!{\rm End}_{R}(P)\ot {\rm End}_{ D_g\ot  D_{gh}\ot D_{ghl}}\!(\! D_g\ot  D_{gh}\ot  J_{ghl})\!\ot {\rm End}_{R}(P^*).$$
Finally, since  the $R$-algebra isomorphisms ${\rm End}_{R}(P)\cong R\cong {\rm End}_{R}(P^*),$ send  ${\rm id}_{P} $ and ${\rm id}_{P^*}$  to $1_R,$
we obtain that  $\om'$  coincides with $\om.$ This shows  that  $\f_6$  is well defined. \cua

\vd

\begin{teo} $\f_6\colon H^1(G,\af^*, {\bf PicS}(R))\to H^3(G,\af, R)$ is a group homomorphism.
\end{teo}
\p   Let $f,f'\in Z^1(G,\af^*, {\bf PicS}(R)).$  Write  $\Phi _f (g)=[J_g]$ and $\Phi  _{f'}(g)=[J'_g].$ Notice that
$ f(g)=f(g)[D_g]=f(g)\Phi_0(g)\Phi_0(g\m)=[J_g\ot {}_{g\m}(D_g)_{I}].
$ Then $\Phi_{f f'} (g)=[J_g\ot {}_{g\m}(D_g)_{I}\ot  J'_g],$ for all $g\in G, $ and
 there are $R$-$R$-bimodule isomorphisms
$$\{F_{g,h}\colon T_g\ot T_h\to  D_g\ot T_{gh}\}_{g,h\in G},$$ where $T_g=J_g\ot{}_{g\m}(D_g)_{I}\ot J'_g,$ $g\in G.$ We shall
make a specific choice of the $F_{g,h}.$  Notice first that
\begin{align*}T_g\ot T_h& \stackrel{(\ref{jgdgi})}{\cong}
 (J_g\ot {}_{g\m}(D_g)_{I}\ot  J'_g \ot D_{g\m})\ot[\underbrace{(J_h\ot {}_{h\m}(D_h)_{I})}_{\in {\bf PicS}(R)}\ot J'_h] \\
&\cong (J_g\ot {}_{g\m}(D_g)_{I}\ot  J'_g \ot D_{g\m})\ot [(J_h\ot {}_{h\m}(D_h)_{I})\ot D_{g\m}\ot  J'_h]. \end{align*}
Moreover,  ${}_{g\m}(D_g)_{I}
\ot {}_{g}(D_{g\m})_{I} \cong D_{g\m}$  by (\ref{dg}), and  we see that  $T_g\ot T_h$ is isomorphic to
\begin{align*} &[(J_g\ot {}_{g\m}(D_g)_{I}\ot \underbrace{(J'_g\ot {}_{g\m}(D_g)_{I}}_{\in\, {\bf PicS}(R)})]
		\ot [\underbrace{{}_{g}(D_{g\m})_{I}\ot  (J_h\ot {}_{h\m}(D_h)_I)\ot _{g\m}(D_g)_{I}}_{\in \,{\bf PicS}(R)}]
&\!\!\!\ot _{g}(D_{g\m})_{I}\ot J'_h\end{align*} as $R$-$R$-bimodules.
Furthermore, since the elements in ${\bf PicS}(R)$ commute, there are $R$-$R$-bimodule isomorphisms
\begin{align*}
T_g\ot T_h&\cong ( J_g\ot {}_{g\m}(D_g)_{I})\ot [{}_{g}(D_{g\m})_{I}\ot (J_h\ot {}_{h\m}(D_h)_I)\ot {}_{g\m}(D_g)_{I}]\ot\\
&\,\,\,\,\,\,\,\,( J'_g\ot{}_{g\m}(D_g)_{I})\ot {}_{g}(D_{g\m})_{I}\ot J'_h\\
&\cong  J_g\ot [{}_{g\m}(D_g)_{I}\ot {}_{g}(D_{g\m})_{I}]\ot [J_h\ot {}_{h\m}(D_h)_I \ot {}_{g\m}(D_g)_{I}]\ot  J'_g \\
&\,\,\,\,\,\,\,\, \ot [{}_{g\m}(D_g)_{I}
\ot {}_{g}(D_{g\m})_{I}]\ot  J'_h\\
&\stackrel{(\ref{dg})}\cong ( J_g\ot  D_{g\m})\ot  J_h\ot{}_{h\m}(D_h)_I \ot {}_{g\m}(D_g)_{I}\ot  (J'_g \ot D_{g\m})\ot  J'_h\\
& \stackrel{(\ref{jgdgi})}{\cong} ( J_g\ot J_h)\ot [{}_{h\m}(D_h)_I\ot {}_{g\m}(D_g)_{I}]
\ot (J'_g \ot J'_h).\end{align*}
Now, applying $\chi _{ g, h  },$  ${\chi }' _{   g, h} $ and Lemma~\ref{isovr}, we get
\begin{align*}
T_{g}\ot T_h&\cong   D_g\ot J_{gh}\ot {}_{h\m}(D_h)_I\ot \underbrace{ ({}_{g\m}(D_g)_{I}\ot   D_g)}\ot J'_{gh}\\
&\cong D_g\ot J_{gh}\ot [{}_{h\m}(D_h)_I\ot {}_{g\m}(D_g)_{I}]\ot J'_{gh}\\
&  {\cong} D_g\ot J_{gh}\ot (D_{h\m}\ot {}_{h\m g\m}(D_{gh})_{I})\ot  J'_{gh}\\
& \cong D_g\ot (J_{gh}\ot  D_{h\m})\ot {}_{h\m g\m}(D_{gh})_{I} \ot  J'_{gh}\\
& \stackrel{(\ref{kappa})}{\cong} D_g\ot  D_g\ot J_{gh}\ot {}_{h\m g\m}(D_{gh})_{I}\ot J'_{gh}\\
&\cong D_g\ot J_{gh}\ot {}_{h\m g\m}(D_{gh})_{I}\ot J'_{gh} = D_g \ot T_{gh},
\end{align*}
and we pick the family $\{F_{g,h}\}_{g,h\in G}$ as the composition of the isomorphisms constructed above. By direct verification we obtain the following:
\begin{claim} The values of $F_{g,h}$ are given by
\begin{equation}\label{F}
F_{g,h}((x_g\ot d_g\ot x'_g)\ot (x_h\ot d_h\ot x'_h))=\chi_{g,h}(x_g\ot x_h)\ot d_g\af_g(d_h1_{g\m})  \chi'_{g,h}(x'_g\ot x'_h),
\end{equation} for any  $(x_g\ot d_g\ot x'_g)\ot (x_h\ot d_h\ot x'_h)\in T_g\ot T_h,\, g,h\in G,$ where $ d_g\af_g(d_h1_{g\m})
\chi'_{g,h}(x'_g\ot x'_h)$ is considered in ${}_{h\m g\m}(D_{gh})_{I}\ot J'_{gh},$   and  is given by
$\sum\limits_{i} d_g\af_g(d_h1_{g\m})e'_{g,i}\ot v'_{gh,i},$ in which ${  \chi ' _{g,h}  }  (x'_g\ot x'_h)=\sum_{i}e'_{g,i}\ot v'_{gh,i}.$
\end{claim}

We shall also need the next.
\begin{claim} The inverse of $F_{g,h}$ is given by \begin{equation}\label{Finv}
d_g\ot x_{gh}\ot d_{gh}\ot x'_{gh}\mt\sum_{i,j}     (y_{g,i}\ot 1_{g}\ot  y'_{g,j})\ot ( z_{h,i}\ot \af_{g\m}(d_{gh}1_g) \ot   z'_{h,j}), \,\,g,h\in G,
\end{equation}
where $\chi_{g,h}(\sum_iy_{g,i}\ot z_{h,i})=d_g\ot x_{gh}$ and $\chi'_{g,h}(\sum_jy'_{g,j}\ot z'_{h,j})=1_g\ot x'_{gh},$ $g,h\in G.$
\end{claim}

Let $V_{g,h}$ be the map defined by \eqref{Finv}. Then,
\begin{align*}
d_g\ot x_{gh}\ot d_{gh}\ot x'_{gh} &\stackrel{V_{g,h}}\mt \sum_{i,j}  (y_{g,i}\ot 1_{g}\ot  y'_{g,j})\ot ( z_{h,i}\ot \af_{g\m}(d_{gh}1_g) \ot   z'_{h,j})\\
&\stackrel{F_{g,h}}\mt \sum_{i,j} [ \chi_{g,h}(y_{g,i}\ot z_{h,i})\ot d_{gh}1_g  \chi'_{g,h}(y'_{g,j}\ot  z'_{h,j})]\\
&=(\sum_{i}\chi_{g,h}(y_{g,i}\ot z_{h,i}))\ot d_{gh}1_g \sum_{j} \chi'_{g,h}(y'_{g,j}\ot  z'_{h,j})\\
&=d_g\ot x_{gh}\ot d_{gh}1_g\ot x'_{gh}\\
&=d_g\ot x_{gh}\ot 1_{h\m}\bullet d_{gh}\ot x'_{gh}\\
&=d_g\ot 1_gx_{gh}\ot  d_{gh}\ot x'_{gh}\\
&=d_g\ot x_{gh}\ot  d_{gh}\ot x'_{gh},
\end{align*} and since $F_{g,h}$ is invertible, we conclude that $V_{g,h} =F\m_{g,h}$ for all $g,h\in G.$

\vd

The automorphism $\tilde\om_F$ induced by the family $\{F_{g,h}\}_{g,h\in G}$ is
$\tilde\om_F(g,h,l)=({\rm id}_{D_g}\ot F_{gh,l})\circ(F_{g,h}\ot{\rm id}_l)\circ({\rm id}_g\ot F_{h,l})\m\circ(\kappa_{g,h}\ot
{\rm id }_{hl})\m\circ({\rm id}_{D_{gh}}\ot F_{g,hl})\m\circ({\tau}_{gh,g}\ot {\rm id}_{ghl})\m.$ Hence, for any
$u_{ghl}= 1_{g}\ot1_{gh}\ot x_{ghl}\ot d_{ghl} \ot x'_{ghl}\in {\rm dom}\,\tilde\om_F(g,h,l),$ there is a unique
$\om_F(g,h,l)\in$\linebreak$ \U(D_gD_{gh}D_{ghl}),$$\, g,h,l\in G$, such that  $\tilde\om_F(g,h,l) ( u_{ghl} ) =\om_F(g,h,l)u_{ghl}$.

\vd

\begin{claim} $\om_F=\om\om',$ or equivalently $\tilde\om_F(g,h,l) (u_{ghl}) = \om\om'(g,h,l)u_{ghl},$ for all $g,h,l\in G.$ \end{claim}

First,  we calculate the value $\tilde\om(g,h,l)(1_g\ot1_{gh}\ot x_{ghl}),$ where $x_{ghl}\in J_{ghl}.$ For this we denote
$$\chi\m_{g,hl}(1_g\ot x_{ghl})=\sum_iy_{g,i}\ot z_{hl,i},\,\,\,\,\,\,\,\,\chi\m_{h,l}(1_h\ot z_{hl,i}) =\sum_m u_{h,m}^{(i)}\ot v_{l,m}^{(i)}$$
and
$$\chi_{g,h}(y_{g,i}\ot u_{h,m}^{(i)})=\sum\limits_ks_{g,k}^{(i,m)}\ot t_{gh,k}^{(i,m)},\,\,\,\,\,\,\,\, \chi_{gh,l}
( t_{gh,k}^{(i,m)}\ot v_{l,m}^{(i)} )= \sum_p c_{gh,p}^{(i,m,k)}\ot e_{ghl,p} ^{(i,m,k)}.$$
Then,
\begin{align*}
1_g\ot 1_{gh}\ot x_{ghl}&\mt 1_{gh}\ot 1_{g}\ot x_{ghl}
\mt 1_{gh}\ot \chi\m_{g,hl}(1_g\ot x_{ghl})\\
&=\sum\limits_i (1_{gh}\ot y_{g,i}\ot z_{hl,i})
\mt \sum\limits_i (y_{g,i} \ot 1_{h} \ot z_{hl,i})\\
&\mt \sum\limits_i [y_{g,i} \ot \chi\m_{h,l}(1_{h} \ot z_{hl,i})]
= \sum\limits_{i,m} [y_{g,i} \ot(u_{h,m}^{(i)}\ot v_{l,m}^{(i)})]\end{align*}
\begin{align*}
&\mt \sum\limits_{i,m} [\chi _{g,h}(y_{g,i} \ot u_{h,m}^{(i)})\ot v_{l,m}^{(i)})]
= \sum\limits_{i,m,k} [(s_{g,k}^{(i,m)}\ot t_{gh,k}^{(i,m)})\ot v_{l,m}^{(i)})]\\
&\mt \sum\limits_{i,m,k} [s_{g,k}^{(i,m)}\ot \chi _{gh,l}( t_{gh,k}^{(i,m)}\ot v_{l,m}^{(i)})]\\
&= \sum\limits_{i,m,k,p} [s_{g,k}^{(i,m)}\ot (c_{gh,p}^{(i,m,k)}\ot  { e_{ghl,p}^{(i,m,k)}  }    )]\\
&=  1_g\ot 1_{gh}\ot \sum\limits_{i,m,k,p} (s_{g,k}^{(i,m)}c_{gh,p}^{(i,m,k)}   { e_{ghl,p}^{(i,m,k)} } ).
\end{align*}

Since  $\tilde\om(g,h,l)(1_g\ot1_{gh}\ot x_{ghl})=\om(g,h,l)(1_g\ot1_{gh}\ot x_{ghl}),$ the uniqueness of  $\om(g,h,l)$ implies
\begin{equation}  \label{omega} 1_g\ot 1_{gh}\ot \sum\limits_{i,m,k,p} (s_{g,k}^{(i,m)}c_{gh,p}^{(i,m,k)}  { e_{ghl,p}^{(i,m,k)} }   )=1_g\ot1_{gh}\ot \om(g,h,l)x_{ghl},\end{equation}
for all $g,h,l \in G.$ Analogously, denoting
$$\chi'\m_{g,hl}(1_g\ot {x'}_{ghl})=\sum\limits_j{y'}_{g,j}\ot z'_{hl,j},\,\,\,\,\,\,\,\,\chi'\m_{h,l}(1_h\ot z'_{hl,j}) =
\sum\limits_n {u'}_{h,n}^{(j)}\ot {v'}_{l,n}^{(j)}$$ and
$$\chi'_{g,h}({y'}_{g,j}\ot {u'}_{h,n}^{(j)})=\sum\limits_{k'}{s'}_{g,{k'}}^{(j,n)}\ot {t'}_{gh,{k'}}^{(j,n)},\,\,\,\,\,\,\,\,
\chi'_{gh,l}({ t'}_{gh,{k'}}^{(j,n)}\ot {v'}_{l,n}^{(j)} )= \sum_{p'} {c'}_{gh,p'}^{(j,n,{k'})}\ot {e'}_{ghl,{p'}} ^{(j,n,{k'})},$$
we obtain
 $ 1_g\ot 1_{gh}\ot\sum\limits_{j,n,{k'},{p'}}    {  {s'}_{g, {k'} }^{(j,n)} }   {c'}_{gh,{p'}}^{(j,n,{k'})} {e'}_{ghl,{p'}}^{(j,n,k')}= 1_g\ot 1_{gh}\ot \om'(g,h,l)x'_{ghl},$ which implies

\begin{equation}  \label{omegal} \sum\limits_{j,n,{k'},{p'}} {s'}_{g,k'}^{(j,n)}{c'}_{gh,{p'}}^{(j,n,{k'})} {e'}_{ghl,{p'}} ^{(j,n,k')}= \om'(g,h,l)x'_{ghl}. \end{equation}
Now we use (\ref{F}), (\ref{Finv}) and  diagram (\ref{diag1}) to calculate $\tilde\om_F(g,h,l)u_{ghl}.$ We have that
\begin{align*}
u_{ghl}&\mt 1_{gh}\ot1_{g}\ot x_{ghl}\ot d_{ghl} \ot x'_{ghl}
\mt 1_{gh}\ot F\m_{g,hl}(1_{g}\ot x_{ghl}\ot d_{ghl} \ot x'_{ghl})\\
&=\sum_{i,j} 1_{gh}\ot(y_{g,i}\ot 1_g\ot y'_{g,j})\ot(z_{hl,i}\ot\af_{g\m}(d_{ghl}1_g)\ot z'_{hl,j})\\
&\mt\sum_{i,j} \kappa\m_{g,h}(1_{gh}\ot(y_{g,i}\ot 1_g\ot y'_{g,j}))\ot(z_{hl,i}\ot\af_{g\m}(d_{ghl}1_g)\ot z'_{hl,j}) \\
&= \sum_{i,j} \underbrace{ y_{g,i}\ot 1_g\ot y'_{g,j} }_{\in T_g}\ot 1_{h} 1_{g\m}\ot(z_{hl,i}\ot\af_{g\m}(d_{ghl}1_g)\ot z'_{hl,j})\\
&= \sum_{i,j} y_{g,i}\ot 1_g\ot y'_{g,j}\ot 1_{h} \ot(z_{hl,i}\ot\af_{g\m}(d_{ghl}1_g)\ot z'_{hl,j})\\
&\mt\sum_{i,j} (y_{g,i}\ot 1_g\ot y'_{g,j})\ot F\m_{h,l}[1_{h} \ot (z_{hl,i} \ot  \af_{g\m}(d_{ghl}1_g)\ot z'_{hl,j}) ]\\
&= \sum_{i,j,m,n}(y_{g,i}\ot 1_g\ot y'_{g,j})\ot [(u_{h,m}^{(i)}\ot 1_h\ot {u'}_{h,n}^{(j)})\ot(v_{l,m}^{(i)}\ot\af_{h\m}(\af_{g\m}(d_{ghl}1_g)1_h)\ot {v'}_{l,n}^{(j)})]\end{align*}
\begin{align*}
&\mt\sum_{i,j\atop m,n}F_{g,h}[(y_{g,i}\ot 1_g\ot y'_{g,j})\ot (u_{h,m}^{(i)}\ot 1_h\ot {u'}_{h,n}^{(j)})]\ot(v_{l,m}^{(i)}\ot\af_{h\m}(\af_{g\m}(d_{ghl}1_g)1_h)\ot {v'}_{l,n}^{(j)})]\\
&= \sum_{i,j\atop k,n} \chi _{g,h}(y_{g,i}\ot u_{h,m}^{(i)})\ot \underbrace{1_g} 1_{gh}   \chi '_{g,h}(y'_{g,j}\ot {u'}_{h,n}^{(j)}) \ot (v_{l,m}^{(i)}\ot\af_{h\m}(\af_{g\m}(d_{ghl}1_g)1_h)\ot {v'}_{l,n}^{(j)})\\
&= \sum_{i,j\atop k,n} \chi _{g,h}(y_{g,i}\ot u_{h,m}^{(i)})\ot 1_{gh}   \chi '_{g,h}(y'_{g,j}\ot {u'}_{h,n}^{(j)}) \ot (v_{l,m}^{(i)}\ot\af_{h\m}(\af_{g\m}(d_{ghl}1_g)1_h)\ot {v'}_{l,n}^{(j)})\\
&=  \sum_{i,j\atop m,n} \chi _{g,h}(y_{g,i}\ot u_{h,m}^{(i)})\ot  1_{gh} \chi '_{g,h}(y'_{g,j}\ot {u'}_{h,n}^{(j)})\ot (v_{l,m}^{(i)}\ot\underbrace{\af_{ (gh)\m}(d_{ghl}1_{gh})}_{\in\,\,{}_{l\m}(D_l)_{I} }1_{h\m}\ot  {v'}_{l,n}^{(j)})\\
&= \sum_{i,j\atop m,n} \chi_{g,h}(y_{g,i}\ot u_{h,m}^{(i)})\ot  1_{gh} \chi'_{g,h}(y'_{g,j}\ot {u'}_{h,n}^{(j)})\ot (\underbrace{v_{l,m}^{(i)}}_{\in\,\,J_l }1_{l\m h\m}\ot\af_{ (gh)\m}(d_{ghl}1_{gh})\ot {v'}_{l,n}^{(j)})\\
&=\sum_{i,j\atop m,n} \chi_{g,h}(y_{g,i}\ot u_{h,m}^{(i)})\ot  1_{gh} \chi'_{g,h}(y'_{g,j}\ot {u'}_{h,n}^{(j)})\ot (v_{l,m}^{(i)}\ot\af_{ (gh)\m}(d_{ghl}1_{gh})\ot {v'}_{l,n}^{(j)})\\
&=\sum_{i,j,k\atop n,m,{k'}} [ (s_{g,k}^{(i,m)}\ot t_{gh,k}^{(i,m)})\ot ( 1_{gh} {s'}_{g,{k'}}^{(j,n)}\ot {t'}_{gh,k'}^{(j,n)})\ot (v_{l,m}^{(i)}\ot\af_{ (gh)\m}(d_{ghl}1_{gh})\ot {v'}_{l,n}^{(j)})]\\
&\mt \sum_{i,j,k\atop n,m,{k'}} s_{g,k}^{(i,m)}\ot F_{gh, l} [( t_{gh,k}^{(i,m)}\ot 1_{gh} {s'}_{g,{k'}}^{(j,n)}\ot {t'}_{gh,k'}^{(j,n)})\ot (v_{l,m}^{(i)}\ot\af_{ (gh)\m}(d_{ghl}1_{gh})\ot {v'}_{l,n}^{(j)})] \\
&= \sum_{i,j,k\atop n,m,{k'}} s_{g,k}^{(i,m)}\ot \chi_{gh,l} (t_{gh,k}^{(i,m)}\ot v_{l,m}^{(i)})\ot  {s'}_{g,{k'}}^{(j,n)}\af_{gh}(\af_{ (gh)\m}(d_{ghl}1_{gh}))\chi'_{gh,l}( {t'}_{gh,k'}^{(j,n)}\ot {v'}_{l,n}^{(j)})\\
&= \sum_{i,j,k\atop n,m,{k'}} s_{g,k}^{(i,m)}\ot \chi_{gh,l} (t_{gh,k}^{(i,m)}\ot v_{l,m}^{(i)})\ot  {s'}_{g,{k'}}^{(j,n)}d_{ghl}\chi'_{gh,l}( {t'}_{gh,k'}^{(j,n)}\ot {v'}_{l,n}^{(j)})\\
& = \sum_{i,m,k} s_{g,k}^{(i,m)}\ot \chi_{gh,l} (t_{gh,k}^{(i,m)}\ot v_{l,m}^{(i)})\ot  d_{ghl}\sum_{j,n ,{k'}}{s'}_{g,{k'}}^{(j,n)}\chi'_{gh,l}( {t'}_{gh,k'}^{(j,n)}\ot {v'}_{l,n}^{(j)})\\
& =  1_g\ot 1_{gh}\ot \sum\limits_{i,m,k,p} (s_{g,k}^{(i,m)}c_{gh,p}^{(i,m,k)} e_{ghl,p}^{(i,m,k)})\ot  d_{ghl}\ot  \sum\limits_{j,n,k',p'} ({s'}_{g,{k'}}^{(j,n)}{c'}_{gh,{p'}}^{(j,n,k')} {e'}_{ghl,{p'}} ^{(j,n,k')})\\
&\stackrel{(\ref{omega},\ref{omegal})}=  1_g\ot 1_{gh}\ot \om(g,h,l)x_{ghl}\ot  d_{ghl}\ot \om'(g,h,l)x'_{ghl}\\
&=\om(g,h,l) \om'(g,h,l) (1_g\ot 1_{gh}\ot x_{ghl}\ot  d_{ghl}\ot x'_{ghl})\\
&=\om(g,h,l) \om'(g,h,l) u_{ghl}.
\end{align*}
Therefore, $\f_6$ is a group homomorphism.
\cua


\end{document}